\documentclass[11pt]{article}   
\usepackage[a4paper,bindingoffset=0.2in,%
            left=1in,right=1in,top=1in,bottom=1in,%
            footskip=.25in]{geometry}    
\usepackage{graphicx}
\usepackage{lineno}
\usepackage{xcolor}
\usepackage{epsfig}
\usepackage[english]{babel}
\usepackage{amsthm,amssymb,amsmath}
 \usepackage{bm}

\newtheorem{theorem}{Theorem}
\newtheorem*{theorem*}{Theorem}

\begin{document}
\begin{center}Hyperbolic trigonometric functions as approximation kernels and their properties I: generalised Fourier transforms
\end{center}

\begin{center} \it Martin Buhmann, Department of Mathematics, Justus-Liebig University of 
   Giessen, Germany, \\ \smallskip
Joaqu\'in J\'odar,
Department of Mathematics, University of Ja\'en, Spain,\\ \smallskip
Miguel L. Rodr\'iguez,
Department of Applied Mathematics, University of Granada, Spain.

\end{center}
\bigskip\bigskip

\today

\par\bigskip\noindent
Abstract: In this paper a new class of radial basis functions based on hyperbolic trigonometric functions will be introduced and studied. We focus on the
 properties of their generalised Fourier transforms with
 asymptotics. Therefore we will compute the expansions of these
 Fourier transforms with an application of the conditions of Strang and
 Fix in order to prove polynomial exactness of
 quasi-interpolants. These quasi-interpolants will be formed with
 special linear combinations of shifts of the new radial functions and
 we will provide explicit expressions for their coefficients. In
 establishing these new radial basis functions we will also use other,
 new classes of shifted thin-plate splines and multiquadrics of \cite{ort2}, \cite{ort3}. There are numerical examples and comparisons of different constructions of quasi-interpolants, in several dimensions, varying the underlying radial basis functions.
\section{Introduction}
For spaces of functions of many variables, one of the most important
challenges is to approximate them by discrete, often finite, sums by
shifts or scaled shifts of kernel functions. The reason why we are
interested in doing that is the use of such approximations for fast
evaluation, efficient combinations with other numerical problems --
for instance integration or computer-aided geometric design -- or
learning with neural networks.

Most of these types of approximations use discrete, semi-discrete or continuous
convolutions with kernel functions. These may be shifted or combined
with distance functions on compact metric spaces of many sorts. The
main point is that we wish to make such approximations available for
any number of variables.

One choice we are particularly interested in and which is one of the
cases of great interest in the mathematical community is the Ansatz
via radial basis functions (RBF). They can be applied for the purpose of
interpolation (see \cite{buh},  for one possible useful reference and a
general overview) or quasi-interpolation (see the recent work \cite{buhjag}). 
In this paper we shall introduce a new set of
kernel functions which we find to be particularly useful because they
combine desirable features of other types of kernel functions. They
are related to tangents hyperbolicus maps but are not identical to
them.

\def\RR{{\mathbb R}}
\def\dd{{\sf d}}
For later reference we note that approximations by radial basis
function kernel functions $\phi:\RR_+\to\RR$ will be written as
\begin{equation}\label{genRBFapprox}
  s(x)=\sum_{\xi\in\Xi}\mu_\xi\phi(\dd(x,\xi)),\qquad x\in\RR^n,
\end{equation}
where $\dd:\RR^n\times\RR^n\to\RR_+$ is a distance function.
$\Xi\subset\RR^n$ is discrete, often consisting of equally spaced
points with respect to the distance function but need not be
finite. Sometimes $\Xi$ has a lattice structure.

Of course, the most frequently used distance function is the Euclidean
norm, but other distance functions that are tailored for special
underlying geometries have been suggested by \cite{BuXu} and they even
render the radial basis functions then used positive definite or
strictly positive definite. These distance functions include the
geodesic distance $\arccos\langle\cdot,\rangle$ on the sphere, or, for
example $$\dd(x,y)=\arccos\left(\langle
x,y\rangle+\sqrt{1-\|x\|^2}\sqrt{1-\|y\|^2}\right)$$ on the unit ball.

On other words, the square interpolation matrices
$$\left(\phi(\dd(\xi,\zeta))\right)_{\xi,\zeta\in\Xi}$$
are positive definite (this is the case when the radial kernel is strictly positive definite, so in particular, non-singular) or positive semi-definite. The latter is so when the radial kernel is positive definite. 

The best-known radial basis functions kernels whose properties we wish
to combine in order to find and analyse a new class include
\begin{itemize}
  \item odd and other non-even (this is highly relevant) powers of the
    distance $\phi(r)=r^\beta$ or sometimes $\phi(r)=r^{2m+1}$,
    \item even powers of distances times logarithms:
      $\phi(r)=r^{2m}\log r$, they are normally called thin-plate
      splines whereas the sheer powers have no own name,
      \item multiquadrics and inverse multiquadrics
        $\phi(r)=\sqrt{r^2+c^2}$ and $\phi(r)=1/\sqrt{r^2+c^2}$ with a
        nonzero constant $c$, and
        \item the new generalised multiquadrics functions and
          thin-plate spline functions (see \cite{ort2}, \cite{ort3}) 
          $$\phi(r)=\left(r^{2\boldsymbol\beta}+c^{2\boldsymbol\beta}\right)^{\boldsymbol\gamma}$$ and
          $$\phi(r)=\left(r^{2\boldsymbol\beta}+c^{2\boldsymbol\beta}\right)^{\boldsymbol\gamma}\log\left(r^{2\boldsymbol\beta}+c^{2\boldsymbol\beta}\right).$$
          They can be used for some (but not all) selections of
          parameters $\boldsymbol\beta$ and $\boldsymbol\gamma$ and for all real $c$. The
          choice of $\boldsymbol\beta=1$ are of course the classical cases.
\end{itemize}
For all of these functions quasi-interpolants can be created which
have the same basic structure as \eqref{genRBFapprox} but use an
intermediate kernel called a quasi-Lagrange function (in a similar use
as the classical Lagrange function except that it does not necessarily
interpolate)
\begin{equation}\label{genQLagrange}
  \psi(x)=\sum_{\xi\in\Xi}\lambda_\xi\phi(\dd(x,\xi)),\qquad x\in\RR^n.
\end{equation}
The coefficients $\lambda_\cdot$ are usually compactly supported with respect to
their index. One then forms the approximation
\begin{equation}\label{genQIapprox}
  s(x)=\sum_{\xi\in\Xi}f(\xi)\psi(\dd(x,\xi)),\qquad x\in\RR^n,
\end{equation}
where $f:\RR^n\to\RR$ is the function to be approximated which needs to be at a minimum continuous for point evaluation.

The most important properties of the quasi-Lagrange functions are
\begin{enumerate}
  \item Polynomial recovery which includes constant polynomials -- the
    quasi-Lagrange functions must form a partition of unity {\it and}
\begin{equation*}
  s(x)=\sum_{\xi\in\Xi}f(\xi)\psi(\dd(x,\xi))\equiv f(x),\qquad x\in\RR^n,
\end{equation*}
for all multivariate polynomials $f$ of some maximal total degree,
\item and combined with this some order of decay
  $$ \psi(r)=O(r^{-\delta})$$
  so as to render the series in the first item well-defined and to
  avail us of some locality in order to be able to establish error
  bounds of some meaningful order.
\end{enumerate}
Both properties can be identified by using the results of computations of generalised
Fourier transforms of the radial kernels and from applications of the
famous Strang and Fix conditions. Roughly speaking, the order of
polynomial exactness comes from the order of the singularities of
those generalised Fourier transforms which in turns depends on the
speed of increase of the radial basis function, and the second
property comes from the smoothness of the quasi-Lagrange function's
Fourier transform which is related to $\phi$'s behaviour at the
origin.

This brings us to our new choices of radial functions which will have,
in the most general expression, the form
\begin{equation}\label{newRBF} \phi(r)=r^\beta\tanh^\alpha r.\end{equation}
We shall also use
\begin{equation}\label{newRBF2} \phi(r)=r^\beta\log r\tanh^\alpha r.\end{equation}
The latter is whenever $\beta$ is an even positive integer. With univariate
applications in mind,  these functions for
the case $\beta=1$, $n=1$ and $\alpha=1$ have been used in \cite{mar}.

We shall only require that $\alpha+\beta>0$. Different cases we shall distinguish
due to the tailored choice of coefficients are when $\beta$ is an integer with separate cases for even and odd numbers.

These functions behave near the origin as a power $r^{\alpha+\beta}$
and towards infinity as a power $r^\beta$. The importance from this
combination comes from the fact explained in the previous paragraph
and in particular from the fact that $r^\beta$ and $r^{\alpha+\beta}$
need not be the same and from the different behaviours of the sheer powers $r^\alpha$ and $r^{\alpha+\beta}$ for approximation when they are compared.

This is because the well-known multiquadric approximations and their
ilk behave asymptotically as a power of $r$ towards infinity but are
infinitely often smooth at the origin. This makes them too inflexible
near zero because as is the case with polynomial splines we do not
wish as a rule to approximate arbitrary approximands with an
infinitely differentiable approximant (hence the success of splines
which are only piecewise smooth). Since $r^{\alpha+\beta}$ is not
  infinitely often differentiable at zero our new classes give these
  approximants exactly those desirable properties of spline-type
  functions. Of course we do not wish to have singularities of the
  radial basis functions in the real domain which is why we require $\alpha+\beta>0$.

  On the other hand, the order of the singularity of the kernel's
  generalised Fourier transform depends on the growth rate of the
  radial basis function. And this singularity, in turn, dictates the
  order of polynomial reproduction of the quasi-interpolant
  \eqref{genQIapprox}.

  Unlike a sheer power (or thin-plate splines)
  this is not fixed together with the smoothness at zero, because, as
  stated earlier $\alpha+\beta$ and $\beta$ are not the same (trivial
  cases being excluded).

  So these new radial basis functions combine the
best of the different kernel functions due to the tangents
hyperbolicus' behaviour at zero and at infinity which is different
asymptotically. This we wish to exploit in this article and it is
evident that the computation of the radial functions generalised
Fourier transform and its behaviour at zero is paramount.

Interestingly enough there is no need to compute explicit expressions
for the generalised Fourier transforms of \eqref{newRBF} (however
mathematically appealing that might be) but what we need are local
expansions in neighbourhoods of the origin to identify the singular
parts. Also, the precise expressions at infinity are not needed, we
merely wish to know the asymptotic orders of decay. These two features
are the properties we shall look for in this paper and we shall
identify them within a great selection of different special cases.

\section{Generalised Fourier transform}\label{sect_particular}

\subsection{Arbitrary dimension}

Let us define $r:= \| {{}}{x} \|$ (${{}}{x}\in\mathbb{R}^n$)
and the radial basis function $g({{}}{x}):=r^m \tanh r$, $m\in
\mathbb{N}$. Then we can prove the following
\begin{theorem} As a radially symmetric function, $g$ has a
  generalised Fourier transform with a short asymptotic expansion at
  zero which has the form
  
  \begin{equation}\label{FirstFT}
  \hat{g}(y)=\hat{u}(y)+ 2^{-m+1}\pi^{\frac{n-1}{2}}  \sum_{j=0}^\infty p_j \|y\|^{2j}, \qquad y\in\mathbb{R}^n,
  \end{equation}
  where  $\hat{u}(y)$ is the generalized Fourier transform of  $r^m$ and which can be found in the appendix.
 The series
on the right-hand side of \eqref{FirstFT} is absolutely convergent in a unit ball around
the origin.
\end{theorem}
\begin{proof}
We begin with rewriting
$$g({{}}{x})=r^m (1-1+\tanh r)=u({{}}{x})-v({{}}{x}),$$
where $u$ is the generalised function $u({{}}{x})=r^m$ and $v$ the $L^1(\mathbb{R}^n)$-function $v({{}}{x})=r^m (1-\tanh r)$.

Therefore, the generalised Fourier transform of $g$, namely $\hat g$, can be split using just for the moment the notation of Jones \cite{jones1982} as follows:
\begin{equation}\label{Fourier_g0}
\hat g({{}}{y})={{}}{\int_{-\infty}^{\infty}} g({{}}{x}) \mathrm{e}^{-{\rm i}
  {{}}{y\cdot x}} \mathrm{d}{{}}{x}={{}}{\int_{-\infty}^{\infty}} u({{}}{x})
\mathrm{e}^{-{\rm i} {{}}{y\cdot x}} \mathrm{d}{{}}{x}-{{}}{\int_{-\infty}^{\infty}}
v({{}}{x}) \mathrm{e}^{-{\rm i} {{}}{y\cdot x}} \mathrm{d}{{}}{x}, \quad y\in\RR^n,
\end{equation}
but we will give up the notation with the non-convergent integrals of
Jones and
write
\begin{equation}\label{Fourier_g}
\hat g(y)=
\hat u(y)
-{{}}{\int_{\RR^n}}
v({{}}{x}) \mathrm{e}^{-{\rm i} {{}}{y\cdot x}} \mathrm{d}{{}}{x}, \qquad y\in\RR^n,
\end{equation}
instead, where only generalised Fourier transforms $\hat g$ and $\hat u
$ are used. The space $\RR^n$ will have to be replaced by the
punctured space $\RR^n\setminus\{0\}$ in due course.

The first expression on the right-hand side is the distributional
Fourier transform of $r^m$, whose value can be obtained from
\cite{jones1982} (see the appendix). On the other hand, since 
$$1-\tanh r=1-\frac{\mathrm{e}^r-\mathrm{e}^{-r}}{\mathrm{e}^r+\mathrm{e}^{-r}}=\frac{2\mathrm{e}^{-r}}{\mathrm{e}^r+\mathrm{e}^{-r}}=\frac{2\mathrm{e}^{-2r}}{1+\mathrm{e}^{-2r}},$$
and taking into account that $v\in L^1(\mathbb{R}^n)$ and \cite[Th. 7.30]{jones1982}, the last integral in (\ref{Fourier_g}) is equal to
\begin{equation}\label{Fourier_f}
{{}}{\int_{\RR^n}} v({{}}{x}) \mathrm{e}^{-{\rm i} {{}}{y\cdot x}} \mathrm{d}{{}}{x}=2^{1+\frac{n}{2}}\pi^{n/2} \|y\|^{1-\frac{n}{2}}\int_0^\infty r^{\frac{n}{2}}r^m\dfrac{\mathrm{e}^{-2r}}{1+\mathrm{e}^{-2r}} J_{\frac{n}{2}-1}(\|y\|r) \mathrm{d}r.
\end{equation}
Using $s$ for $\|y\|$, $y\in\RR^n\setminus\{0\}$, and continuing using
$r$ for $\|x\|$, $x\in\RR^n$, we have that
\begin{align*}
\int_0^\infty r^{m+{\frac{n}{2}}} \dfrac{\mathrm{e}^{-2r}}{1+ \mathrm{e}^{-2r}} J_{\frac{n}{2}-1}(sr) \mathrm{d}r
&=\int_0^\infty r^{m+{\frac{n}{2}}} \mathrm{e}^{-2r}\sum_{k=0}^{\infty}(-1)^k\mathrm{e}^{-2rk} J_{\frac{n}{2}-1}(sr) \mathrm{d}r\\ 
&=\int_0^\infty r^{m+{\frac{n}{2}}} \mathrm{e}^{-2r}\sum_{k=0}^{\infty}(-1)^k\mathrm{e}^{-2rk}
\left(\dfrac{s r}{2}\right)^{\frac{n}{2}-1}\\
&\sum_{j=0}^\infty \dfrac{(-1)^j}{j! \Gamma(\frac{n}{2}+j)}\left(\dfrac{s r}{2}\right)^{2j}\mathrm{d}r,\\
\end{align*}
where the expression  8.440 of \cite{grad} has been used for the expansion of $J_{\frac{n}{2}-1}(sr).$

Now, 
\begin{align*}
&\int_0^\infty r^{m+{\frac{n}{2}}} \mathrm{e}^{-2r}\sum_{k=0}^{\infty}(-1)^k\mathrm{e}^{-2rk}
\left(\dfrac{s r}{2}\right)^{\frac{n}{2}-1}\sum_{j=0}^\infty \dfrac{(-1)^j}{j! \Gamma(\frac{n}{2}+j)}\left(\dfrac{s r}{2}\right)^{2j}\mathrm{d}r\\
&=\sum_{j=0}^{\infty} \sum_{k=0}^\infty(-1)^k \dfrac{(-1)^j}{j! \Gamma(\frac{n}{2}+j)}\left(\dfrac{s}{2}\right)^{2j+\frac{n}{2}-1} \int_0^\infty r^{m+ 2j+n-1} \mathrm{e}^{-2r(k+1)}
\mathrm{d}r\\
&= \sum_{j=0}^\infty \sum_{k=1}^{\infty} \dfrac{(-1)^{k+j-1}}{j! \Gamma(\frac{n}{2}+j)}\left(\dfrac{s}{2}\right)^{2j+\frac{n}{2}-1} \int_0^\infty r^{m+ 2j+n-1} \mathrm{e}^{-2rk}
\mathrm{d}r\\
&= \sum_{j=0}^\infty\sum_{k=1}^{\infty}\dfrac{(-1)^{k+j-1}}{j! \Gamma(\frac{n}{2}+j)}\left(\dfrac{s}{2}\right)^{2j+\frac{n}{2}-1} 
\dfrac{\Gamma(m+ 2j+n)}{(2k)^{m+ 2j+n}}\\
&= \sum_{j=0}^\infty\sum_{k=1}^{\infty}\dfrac{(-1)^{k+j-1}}{j! 2^{m+4j+\frac{3n}{2}-1} k^{m+ 2j+n}}
\dfrac{\Gamma(m+ 2j+n)}{\Gamma(\frac{n}{2}+j)} s^{2j+\frac{n}{2}-1},
\end{align*}
where the expression  3.381 (4) of \cite{grad}  has been used in the last integral of the display.
Since
\begin{align*}
\dfrac{\Gamma(m+ 2j+n)}{\Gamma(\frac{n}{2}+j)} &= \dfrac{\Gamma(2j+n) \prod_{i=0}^{m-1}(n+2j +i)}{\Gamma(\frac{n}{2}+j)}\\
&=\dfrac{2^{2j+n-1} \Gamma(j+\frac{n}{2}) \Gamma(j+\frac{n}{2}+ \frac{1}{2})}{\Gamma(\frac{n}{2}+j)\sqrt{\pi}}\prod_{i=0}^{m-1}(n+2j +i)\\
&=\dfrac{2^{2j+n-1}\Gamma(j+\frac{n}{2}+ \frac{1}{2})}{\sqrt{\pi}}\prod_{i=0}^{m-1}(n+2j +i),
\end{align*}
(where the expression 8.335 (1) \cite{grad} has been used in the second line of the display) we carry on with
\begin{align*}
& \sum_{j=0}^\infty\sum_{k=1}^{\infty}\dfrac{(-1)^{k+j-1}}{j! 2^{m+4j+\frac{3n}{2}-1} k^{m+ 2j+n}}
\dfrac{\Gamma(m+ 2j+n)}{\Gamma(\frac{n}{2}+j)} s^{2j+\frac{n}{2}-1}\\
&=
\sum_{j=0}^\infty\sum_{k=1}^{\infty} \dfrac{(-1)^{k+ j-1}}{j! k^{m+ 2j+n} 2^{m+2j+\frac{n}{2}}}
  \prod_{i=0}^{m-1}(n+2j +i)\dfrac{\Gamma(j+\frac{n+1}{2})}{\sqrt{\pi}} s^{2j+\frac{n}{2}-1} \\
  &=
  \dfrac{-1}{2^{m+\frac{n}{2}}\sqrt{\pi}}\sum_{j=0}^\infty \left(\sum_{k=1}^{\infty} \dfrac{(-1)^{k}}{k^{m+ 2j+n}}\right)
 \dfrac{ (-1)^{j}\displaystyle\prod_{i=0}^{m-1}(n+2j +i) \Gamma\left(j+\frac{n+1}{2}\right)}{2^{2j} j!} s^{2j+\frac{n}{2}-1}.
\end{align*}

If we define $$p_j= \dfrac{(-1)^j }{2^{2j} j!}\left(\sum_{k=1}^{\infty} \dfrac{(-1)^{k}}{k^{m+ 2j+n}}\right)\prod_{i=0}^{m-1}(n+2j +i) \Gamma(j+\frac{n+1}{2}),\quad  j\geq 0$$
the expression (\ref{Fourier_f}) can be written as 
$$
{{}}\int_{\RR^n} v({{}}{x}) \mathrm{e}^{-{\rm i} {{}}{y\cdot x}} \mathrm{d}{{}}{x}=-2^{-m+1}\pi^{\frac{n-1}{2}} s^{1-\frac{n}{2}}  \sum_{j=0}^\infty p_j s^{2j+\frac{n}{2}-1}.
$$
Therefore:
$$
\hat g(y)=\hat u(y)+
2^{-m+1}\pi^{\frac{n-1}{2}}  \sum_{j=0}^\infty p_j s^{2j}
$$ 
which will work for $|s|<1$ because then the last series converges;
otherwise the power series may diverge.
\end{proof}

\subsection{Odd dimension}
We have just seen that the above study does not provide an explicit expression for $|s|\geq1.$  In this subsection we develop a valid expression for all $s\in\mathbb{R}.$ For this end,  we restrict the framework to odd dimensions, i.e., to $\mathbb{R}^n$ with $n$ an odd natural number.  The advantage of this is that  we will deal with a finite expression of the Bessel functions which help us in order to study the convergence of the resulting series. We will
continue with the other cases, e.g., non-integer degree singularities at zero later-on.

We start with
\begin{align*}
\int_0^\infty r^{m+{\frac{n}{2}}} \dfrac{\mathrm{e}^{-2r}}{1+ \mathrm{e}^{-2r}} J_{\frac{n}{2}-1}(sr) \mathrm{d}r
&=\int_0^\infty r^{m+{\frac{n}{2}}}\sum_{k=0}^{\infty}(-1)^k \mathrm{e}^{-2r(k+1)} J_{\frac{n}{2}-1}(sr) \mathrm{d}r\\
& =\sum_{k=0}^{\infty}(-1)^k\int_0^\infty r^{m+{\frac{n}{2}}}\mathrm{e}^{-2r(k+1)} J_{\frac{n}{2}-1}(sr) \mathrm{d}r.
\end{align*}
If $n=1$ then the formula 8.464 (2) \cite{grad} is
$$J_{\frac{n}{2}-1}(sr)=J_{-\frac{1}{2}}(sr)=\sqrt{\frac{2}{\pi s r}} \cos (sr),$$
and then
\begin{small}
\begin{align*}
\int_0^\infty r^{m+{\frac{n}{2}}} \dfrac{\mathrm{e}^{-2r}}{1+ \mathrm{e}^{-2r}} J_{\frac{n}{2}-1}(sr) \mathrm{d}r
& =\sum_{k=0}^{\infty}(-1)^k \sqrt{\frac{2}{\pi s}}\int_0^\infty r^m \mathrm{e}^{-2r(k+1)} \cos (sr) \mathrm{d}r\\
=& \sum_{k=0}^{\infty}(-1)^k \sqrt{\frac{2}{\pi s}}
\dfrac{\Gamma(m+1)}{(4(k+1)^2+ s^2)^{\frac{m+1}{2}}}\times\\
&\times\cos\left((m+1)\arctan \left(\dfrac{s}{2(k+1)}\right)\right),
\end{align*}
\end{small}
where the expression 3.944 (6) of \cite{grad} has been used into the last integral.

 On the other hand, if $n\geq 3$, 
 by applying 8.461 (1) of \cite{grad}, 
\begin{align*}
J_{\frac{n}{2}-1}(sr)=&
\sqrt{\dfrac{2}{\pi sr}} \left(\sin(sr- \frac{\pi(n-3)}{4}) \displaystyle\sum_{\ell=0}^{[{\frac{n-3}{4}}]}\dfrac{(-1)^\ell \left(\frac{n-3}{2}+2\ell\right)!}{(2\ell)! \left(\frac{n-3}{2}-2\ell\right)!}(2sr)^{-2\ell}\right.
\\ & +\left.\cos\left(sr- \frac{\pi(n-3)}{4}\right)
\displaystyle\sum_{\ell=0}^{[{\frac{n-5}{4}}]}\dfrac{(-1)^\ell\left(\frac{n-3}{2}+2\ell+1\right)!}{(2\ell+1)! \left(\frac{n-3}{2}-2\ell-1\right)!}(2sr)^{-(2\ell+1)}\right).  
\end{align*}
The square
brackets denote the Gau\ss bracket, i.e., the largest integer at most
its argument.

Denoting by $(a\mod b)$ the remainder when dividing $a$ by $b$ by the division algorithm, the last expression can be simplified depending on the value of the nonnegative integer $\frac{n-3}{2}$ as follows:
\begin{itemize}
\item If $(\frac{n-3}{2}\mod 4)$ is even, then
$$
J_{\frac{n}{2}-1}(sr)
=\sum_{\ell =0}^{[\frac{n-3}{4}]}C_\ell \; s^{-2\ell -\frac{1}{2}} \; r^{-2\ell-\frac{1}{2}} \sin (sr)
+\sum_{\ell =0}^{[ \frac{n-5}{4} ]} D_\ell \; s^{-2\ell-\frac{3}{2}} \; r^{-2\ell-\frac{3}{2}} \; \cos (sr),
$$
where
$$C_\ell =(-1)^{\frac{(\frac{n-3}{2}\mod 4)}{2}} \, 2^{-2\ell+\frac{1}{2}}\pi^{-1/2}\frac{(-1)^\ell (\frac{n-3}{2}+2\ell)!}{(2\ell)!(\frac{n-3}{2}-2\ell)!},$$
and
$$
D_\ell = (-1)^{\frac{(\frac{n-3}{2}\mod 4)}{2}}\, 2^{-2\ell-\frac{1}{2}} \pi^{-1/2} \frac{(-1)^\ell (\frac{n-3}{2}+2\ell+1)!}{(2\ell+1)!\left(\frac{n-3}{2}-2\ell-1\right)!},
$$
and
\begin{small}
\begin{align*}
&\int_0^\infty r^{m+{\frac{n}{2}}} \dfrac{\mathrm{e}^{-2r}}{1+ \mathrm{e}^{-2r}} J_{\frac{n}{2}-1}(sr) \mathrm{d}r\\
& =\sum_{k=0}^{\infty}(-1)^k \sum_{\ell =0}^{[\frac{n-3}{4}]}C_\ell \; s^{-2\ell -\frac{1}{2}} \int_0^\infty r^{m+\frac{n}{2}-2\ell-\frac{1}{2}}  \mathrm{e}^{-2r(k+1)} \sin (sr) \mathrm{d}r\\
&+\sum_{k=0}^{\infty}(-1)^k \sum_{\ell =0}^{[ \frac{n-5}{4} ]} D_\ell \; s^{-2\ell-\frac{3}{2}} \; \int_0^\infty r^{m+\frac{n}{2}-2\ell-\frac{3}{2}} \mathrm{e}^{-2r(k+1)}  \cos (sr)\mathrm{d}r\\
&=\sum_{k=0}^{\infty}(-1)^k \sum_{\ell
    =0}^{[\frac{n-3}{4}]}C_\ell \; s^{-2\ell -\frac{1}{2}}
  \dfrac{\Gamma( m+\frac{n}{2}-2\ell+\frac{1}{2} )}{(4(k+1)^2+
    s^2)^{\frac{2m+n+1}{4}-\ell}}\times\\
  &\times\sin\left(\left({\frac{2m+n+1}{2}-2\ell}\right)\arctan
  \left(\dfrac{s}{2(k+1)}\right)\right)\\ 
&+\sum_{k=0}^{\infty}(-1)^k \sum_{\ell =0}^{[ \frac{n-5}{4}
    ]} D_\ell \; s^{-2\ell-\frac{3}{2}} \dfrac{\Gamma(
    m+\frac{n}{2}-2\ell-\frac{1}{2} )}{(4(k+1)^2+
    s^2)^{\frac{2m+n-1}{4}-\ell}}\times\\
  &\times\cos\left(\left({\frac{2m+n-1}{2}-2\ell}\right)\arctan \left(\dfrac{s}{2(k+1)}\right)\right)
\end{align*}
\end{small}

\item If $(\frac{n-3}{2}\mod 4)$ is odd, then
\begin{equation*}
J_{\frac{n}{2}-1}(sr)
=\sum_{\ell =0}^{[\frac{n-3}{4}]}E_\ell \; s^{-2\ell -\frac{1}{2}} \; r^{-2\ell-\frac{1}{2}} \cos (sr)
+\sum_{\ell =0}^{[ \frac{n-5}{4} ]} F_\ell \; s^{-2\ell-\frac{3}{2}} \; r^{-2\ell-\frac{3}{2}} \; \sin (sr),
\end{equation*}
where
$$
E_\ell=(-1)^{\frac{(\frac{n-3}{2}\mod 4)+1}{2}}\,2^{-2\ell+\frac{1}{2}}\pi^{-1/2}\frac{(-1)^\ell (\frac{n-3}{2}+2\ell)!}{(2\ell)!(\frac{n-3}{2}-2\ell)!},
$$
and
$$F_\ell=(-1)^{\frac{(\frac{n-3}{2}\mod 4)-1}{2}}\, 2^{-2\ell-\frac{1}{2}} \pi^{-1/2} \frac{(-1)^\ell (\frac{n-3}{2}+2\ell+1)!}{(2\ell+1)!\left(\frac{n-3}{2}-2\ell-1\right)!},
$$

and
\begin{small}
\begin{align*}
&\int_0^\infty r^{m+\frac{n}{2}} \dfrac{\mathrm{e}^{-2r}}{1+ \mathrm{e}^{-2r}} J_{\frac{n}{2}-1}(sr) \mathrm{d}r\\
& =\sum_{k=0}^{\infty}(-1)^k \sum_{\ell =0}^{[\frac{n-3}{4}]}E_\ell \; s^{-2\ell -\frac{1}{2}}\int_0^\infty r^{m+\frac{n}{2}-2\ell -\frac{1}{2}} \mathrm{e}^{-2r(k+1)} \cos (sr) \mathrm{d}r\\
&+\sum_{k=0}^{\infty}(-1)^k \sum_{\ell =0}^{[ \frac{n-5}{4} ]} F_\ell \; s^{-2\ell-\frac{3}{2}}\int_0^\infty r^{m+\frac{n}{2}-2\ell -\frac{3}{2}} \mathrm{e}^{-2r(k+1)} \sin (sr) \mathrm{d}r\\
& =\sum_{k=0}^{\infty}(-1)^k \sum_{\ell
    =0}^{[\frac{n-3}{4}]}E_\ell \; s^{-2\ell -\frac{1}{2}}
  \dfrac{\Gamma(m+\frac{n}{2}-2\ell +\frac{1}{2})}{(4(k+1)^2+
    s^2)^{\frac{2m+n+1}{4}-\ell}}\times\\
  &\times\cos\left((m+\frac{n}{2}-2\ell +\frac{1}{2})\arctan \left(\dfrac{s}{2(k+1)}\right)\right)\\
&+\sum_{k=0}^{\infty}(-1)^k \sum_{\ell =0}^{[ \frac{n-5}{4} ]} F_\ell \; s^{-2\ell-\frac{3}{2}}
\dfrac{\Gamma(m+\frac{n}{2}-2\ell -\frac{1}{2})}{(4(k+1)^2+ s^2)^{\frac{2m+n-1}{4}-\ell}}\times\\&\times\sin\left((m+\frac{n}{2}-2\ell -\frac{1}{2})\arctan \left(\dfrac{s}{2(k+1)}\right)\right)\\
\end{align*}
\end{small}

\end{itemize}

The expressions   3.944 (5) and  3.944 (6)  of \cite{grad} have been used into these integrals similarly to the case $n=1$.

\section{Quasi-interpolation}

\subsection{Classical case}

The study of the polynomial reproduction and order of convergence of
the {\it quasi-interpolant\/} will be based on the fulfillment of the
Strang and Fix conditions which we summarise as follows: 

\begin{theorem}\label{SF}[Strang and Fix conditions]

Let $\psi:\mathbb{R}^n \rightarrow
\mathbb{R}$ be a continuous function such that
\begin{enumerate}
\item there exists a positive $\ell$ such that for some nonnegative integer $M$, when
  $\|{{}}{x}\|\rightarrow \infty$, $|\psi({{}}{x})|={O}(\Vert
  {{}}{x}\Vert^{-n-M-\ell})$, which immediately implies $M$-fold
  differentiability of the Fourier transform,
\item $D^{\alpha} \hat{\psi}(0)=0$, $\forall \alpha \in \mathbb{N}_0^n$, $1\leq \vert \alpha\vert\leq M$, and $\hat{\psi}(0)=1$, where $\vert \alpha \vert =\alpha_1+\cdots+\alpha_n,$
\item $ D^{\alpha}\hat{\psi}(2\pi {{}}{j})=0,\ \forall
  {{}}{j}\in \mathbb{Z}^n\setminus\{0\}$ and  $\forall \alpha \in
  \mathbb{N}_0^n$ with $ \vert\alpha\vert \leq M$.  
\end{enumerate}

Then the quasi-interpolant
\begin{equation}
Q_hf({{}}{x})=\sum_{j\in\mathbb{Z}^n}f(h{{}}{j})\psi({{}}{x}/h-{{}}{j}),\qquad {{}}{x}\in\mathbb{R}^n,
\end{equation}
is well-defined and exact on the linear space $\mathbb{P}_M$ of polynomials of degree less than or equal to $M$. The approximation error can be estimated by
$$\Vert Q_hf-f\Vert_{\infty}=\begin{cases}{O}(h^{M+\ell}), & \text{when } 0<\ell<1,\\
{O}(h^{M+1}\log (1/h)), & \text{when } \ell=1,\\
{O}(h^{M+1}), & \text{when } \ell>1,\\
\end{cases}$$
for $h\rightarrow 0$ and a bounded function $f\in C^{M+1}(\mathbb{R}^n)$ with bounded derivatives.
\end{theorem}

If we consider the function $g({{}}{x})=r^m \tanh r$ from Section
\ref{sect_particular}, we recollect that it can be split as
$g({{}}{x})=u({{}}{x})-v({{}}{x})$, where $u$ is the function
$u({{}}{x})=r^m$ and $v$ the exponentially decaying function when
$r=\|x\|$ tends to $\infty$,  
$$v({{}}{x})=r^m (1-\tanh r)=r^m \frac{2\mathrm{e}^{-2r}}{1+\mathrm{e}^{-2r}}.$$

Let us study now the decay of the generalised function $u$. If
$m\not\in 2\mathbb{N}$, the distributional Fourier transform of $u$ is given by (see the appendix)
\begin{equation}\label{Fourier_t}
 \hat u({{}}{y})=\frac{\Gamma\left(\frac{1}{2}m+\frac{1}{2}n\right)}{\Gamma\left(-\frac{1}{2}m\right)} 2^{m+n}\pi^{n/2}\|y\|^{-m-n}.
\end{equation}

Additionally to this assumption, $m\not\in 2\mathbb{N}$, let us also assume that $m+n\in 2\mathbb{N}$, i.e., $m$ and $n$ are odd natural numbers. In this case, by taking $m+n$ in 
\cite[Theorem 4.8]{buh} 
(subject to a sign change of the function $u$ if needed in order to achieve the
positivity of $\hat u({{}}{y})$ required by (A1) in that theorem,
because of the sign of $\Gamma\left(-\frac{1}{2}m\right)$ in $\hat u (y)$), we
have that there is a multivariate trigonometric polynomial  
\begin{equation}\label{trigPOLY}P({{}}{y})=\sum_{{{}}{k}\in\Delta}
  \mu_{{{}}{k}} \mathrm{e}^{-{\rm{i}}\,{{}}{k}\cdot {{}}{y}},  \end{equation}
being $\Delta$ a finite subset of $\mathbb{Z}^n$, such that the function
$$\psi_u({{}}{x})=\sum_{{{}}{k}\in\Delta}  \mu_{{{}}{k}} u ({{}}{x} - {{}}{k} ),$$
has the asymptotic property
$$\vert \psi_u({{}}{x})\vert={\cal O}(\Vert {{}}{x} \Vert ^{-2n-m})$$
for large $\Vert {{}}{x} \Vert$. According to the proof of that theorem,
$$\hat{\psi_u}({y})=P({y}) \hat u({y})=1+{\cal O}(\Vert {y} \Vert ^{n+m})$$
in a neighbourhood of the origin, and $\hat{\psi_u}$ also satisfies conditions (2) and (3) of Theorem \ref{SF}, with $M=n+m-1$ and $\ell =1$. 

Therefore, by considering
\begin{equation}\label{quasiLF}\psi({{}}{x})=\sum_{{{}}{k}\in\Delta}
  \mu_{{{}}{k}} g ({{}}{x} - {{}}{k}
  )=\psi_u({{}}{x})-\psi_v({{}}{x}),\end{equation} 
where $\psi_v({{}}{x})=\sum_{{{}}{k}\in\Delta}  \mu_{{{}}{k}} v ( {{}}{x} - {{}}{k} )$, and taking into account the exponential decay of $\psi_v({{}}{x})$ for large $\Vert {{}}{x} \Vert$ provided by the one of $v$, we have that
$$\vert \psi({{}}{x})\vert={\cal O}(\Vert {{}}{x} \Vert ^{-2n-m})$$
for large $\Vert {{}}{x} \Vert$. Moreover, taking into account for the properties of the trigonometric polynomial $P$, the distributional Fourier transform of $\psi$,
$$\hat{\psi}({y})=\hat{\psi_u}({y})-\hat{\psi_v}({y})=P({y}) \hat
u({y})-P({y}) \hat v({y}),$$ 
satisfies conditions (2) and (3) of Theorem \ref{SF}, with $M=n+m-1$
and $\ell =1$ (here $\hat v({y})$ is the Fourier transform of the
$L^1(\mathbb{R}^n)$-function $v$). Furthermore, by applying this
theorem, we get the following result.
\begin{theorem} With the construction above, the quasi-interpolant 
$$Q_hf({{}}{x})=\sum_{{{}}{j}\in\mathbb{Z}^n}f(h{{}}{j})\psi({{}}{x}/h-{{}}{j})$$ 
is well-defined and exact on $\mathbb{P}_{n+m-1}$ and the approximation error can be estimated by 
$$\Vert Q_hf-f\Vert_{\infty}={\cal O}(h^{n+m}\log (1/h))$$
for $h\rightarrow 0$ and a bounded function $f\in C^{n+m}(\mathbb{R}^n)$ with bounded derivatives.
\end{theorem}

\subsection{Non-classical case}
If the generalised Fourier transform of the kernel function has an
even order singularity, there is, as we have seen, a finite degree
trigonometric polynomial \eqref{trigPOLY} which has a zero of the same order at the
origin, and we resolve the mentioned singularity by multiplication of
the radial basis function's generalised Fourier transform by that
trigonometric polynomial $P$. Thus the Fourier transform $\hat\psi$ of
the anticipated quasi-Lagrange function $\psi$ is formed, and the
coefficients of the trigonometric polynomial's expansion are the
quasi-Lagrange function's coefficients \eqref{quasiLF}.

Let us call the order of the singularity of $\hat\phi(s)$ as $n+\beta$,
so the generalised Fourier transform has the short asymptotic expansion
\begin{equation}\label{orderSING}
  \hat\phi(s)=c_0 s^{-n-\beta}+c_1 s^{-n-\beta+\delta}+d_0
  s^{-n-\beta+\eta}\log s + \cdots,\qquad s\to0,
\end{equation}
where $n$ still is the space dimension of the ambient space,
$\beta>-n$ is a real parameter, $\delta>0$ marks the next term in the
expansion (usually, $\delta=1$) and $\eta\geq0$ (often $\eta=d_0=0$ and
there is no logarithmic term). The $c_i$ and $d_0$ are real
coefficients, but $c_0$ must not be zero, and as remarked above, $d_0$
is often zero in the prime examples of radial basis functions.

In the extreme situation when there is a logarithmic singularity as a
dominant term, we
allow that either $c_0=0$ or $\beta=-n$ as long as $d_0\neq0$ and $\eta=0$.
\begin{theorem}
  Let $\phi:\RR_+\to\RR$ be a radial basis function whose generalised
  Fourier transform in $n$ dimensions has a short Taylor expansion at zero given by 
  \eqref{orderSING}. Then there is a $2\pi$-periodic function $T$
  such that $T\times\hat\phi(\|\cdot\|):[-\pi,\pi]^n\to\RR$ has a removable
  singularity at zero and its value at zero is one. If
  $\hat\phi(\|\cdot\|)$ restricted to the 
  complement of the unit ball is absolutely integrable, we define
  $$\hat\psi:=T\times\hat\phi(\|\cdot\|):\RR^n\to\RR$$
  which is absolutely integrable and whose inverse Fourier transform
  $\psi$ is continuous. If $\hat\phi$'s singularity is an
  even integer, $T$ may be chosen as a trigonometric polynomial.
  \end{theorem}
  \begin{proof}
The proof of this theorem is in several steps.
We have three cases.
\begin{enumerate}
  \item The first case is the one treated already above,
namely when $n+\beta$ is an even positive integer (and: $d_0=0$ or $\eta>0$). Then the
quasi-Lagrange function's Fourier transform is
$$\hat\psi=P\hat\phi(\|\cdot\|).$$
\item Then there are the other two cases, namely when $n+\beta$ is not an
  even integer -- or not an integer at all -- (and: $d_0=0$  or $\eta>0$), and
  \item when $d_0\neq0=\eta$,
    so that there is a leading order logarithmic term. An example for
    this is the inverse multiquadric function in one dimension where
    $n=1$, $\beta=-1$, $\eta=0$ and $d_0\neq0$. Note that we
    can also allow $d_0\neq0\neq\eta$, and then there is a logarithmic
    singularity but it is not the leading term. A good example for
    that situation is the ubiquitous multiquadric radial function
    where $\beta=\delta=1$, $d_0\neq0$ and $\eta=n+1$. 
\end{enumerate}
We deal with the
cases when $n+\beta$ is not an
even integer (or not an integer at all) and when there is a
logarithmic singularity separately, but we shall consider the two
possibilities when $n+\beta$ is integral but not an
even integer and not even an integer at all at the same time.

If that is so, we form a trigonometric polynomial $P$ as before which
has a zero of even order
$$2\left[\frac{n+\beta}2\right]$$
first in the same ways as we built \eqref{trigPOLY}. 
Note that this number is the same as $n+\beta$ if and only if
$n+\beta$ is an even integer.

Then, for the next step, we follow the construction in the recent paper \cite{jjmm} which
is as follows. Namely, then we create another trigonometric series (but not a
polynomial)
\begin{equation}\label{unevenSING}
  G(y) =\left\{\sqrt{ \sum_{j=1}^n|\sin(y_j)|^2}\right\}^\gamma,\quad
  y=(y_1,y_2,\ldots,y_n),
  \end{equation}
or alternatively (but not equivalently)
\begin{equation}\label{unevenSINGalt}
  G(y) =2^{\frac\gamma2} \left\{\sqrt{\sum_{j=1}^n(1-\cos(y_j))}\right\}^\gamma,\quad
  y=(y_1,y_2,\ldots,y_n),
  \end{equation}
where
$$\gamma=n+\beta-2\left[\frac{n+\beta}2\right]$$
which is always positive under our conditions. It is also strictly
smaller than two.

Then the quasi-Lagrange function's Fourier transform is defined as
\begin{equation}\label{quasiLFnonevenSING}
  \hat\psi=G\times P\times\hat\phi(\|\cdot\|).
  \end{equation}
This follows the construction in \cite[p.~159]{buhjag}.

This periodic function $G$ has an infinite expansion whose coefficients $f_j$, say, are at a
minimum square-summable and whose values are stated in the said book
on these pages in a univariate example. The decay of these
coefficients is also stated in that book as
$$ |f_j|=O\left(\|j\|^{-n-\gamma}\right),\qquad j\neq0, $$
which makes them absolutely summable as long as $\gamma$ is positive
(if it is not, we are back in the classical case).

Of course the coefficients of the quasi-Lagrange function in the real
domain are then formed as the discrete convolution $\underline P*\underline G$ of the coefficient
sequences $\underline P$ and $\underline G$ of $P$ and $G$, respectively.

We note that these convolutions are well defined because $P$ expands
in a trigonometric polynomial (i.e., there are only finitely many
nonzero coefficients) and the infinite expansion of the
$2\pi$-periodic function $G$ has absolutely summable
coefficients. Therefore the resulting convolution also has absolutely
summable coefficients by H\"older's inequality.

If there is a logarithmic singularity, so $\eta=0\neq d_0$, and if $n=1$, \eqref{quasiLFnonevenSING} is
replaced by 
\begin{equation}\label{quasiLFlogSING}
  \hat\psi=G\times H\times P\times\hat\phi(\|\cdot\|).
\end{equation}
where 
$$H(y)= \frac1{d_0}\frac1{-\log\left(\frac12\big|\sin\left(\frac12
  y\right)\big|\right)}$$
  is another $2\pi$-periodic trigonometric function (but not a
  trigonometric polynomial) whose (infinitely many) coefficients
  are at a minimum square summable.

  In fact, it follows from the work
  in \cite[p.~243]{buhjag}, that its coefficients,
  $\underline \mu_j$ say,
  satisfy the bound
  $$|\underline \mu_j|=O\left(\frac1{|j|\log^2|j|}\right),\qquad j\neq -1,0,1.$$
  Therefore the coefficients are actually  absolutely summable.

  Neither $P$, nor $G$, nor $H$ are unique.
  
Of course the coefficients of the quasi-Lagrange function in the real
domain are then formed as the discrete convolution $\underline H*\underline P* \underline G$ of the coefficient
sequences $\underline H$, $\underline P$ and $ \underline G$ of $H$, $P$ and $G$, respectively.

We note once more that these convolutions are well defined and result
in absolutely summable (in particular bounded) coefficients because $P$ expands
in a trigonometric polynomial -- there are only finitely many
nonzero coefficients) -- and the infinite expansions of the
$2\pi$-periodic functions $G$ and $H$ have absolutely summable
coefficients. Therefore the resulting double convolution also has absolutely
summable coefficients by H\"older's inequality.

In all three cases, the generalised Fourier transform $\hat\psi$
satisfies the Strang-and-Fix conditions of order zero, because
$$P(0)=G(0)=H(0)=0,$$
and therefore $\hat\psi(2\pi j)=0$ for all nonzero integral $j$,
but either $P$ or $G$ or $H$ have to be normalised so that
$$\hat\psi(0)=1.$$ The same function $\hat\psi$ is also absolutely
integrable, because $\hat\phi(\|\cdot\|)$ is outside a unit ball about
the origin and because $P$, $G$ and $H$ are continuous periodic functions.
\end{proof}
Examples include the novel radial basis functions of \cite{ort2}
 $$\phi(r)=\left(r^{2\boldsymbol\beta}+c^{2\boldsymbol\beta}\right)^{\boldsymbol\gamma}$$ and of \cite{ort3}
$$\phi(r)=\left(r^{2\boldsymbol\beta}+c^{2\boldsymbol\beta}\right)^{\boldsymbol\gamma}\log\left(r^{2\boldsymbol\beta}+c^{2\boldsymbol\beta}\right)$$
whose generalised Fourier transforms satisfy all our requirements for $n+\beta=n+2\boldsymbol\beta \boldsymbol\gamma$ and decay exponentially, although they are not of one sign for $r>0$.

\section{Extension}

Let us now consider 
$$g(x)= r^{\beta}\left(\tanh r\right)^\alpha, \qquad x\in\mathbb{R}^n, \quad  r:=\Vert x \Vert,$$
 as an extension of the previously considered radial basis functions
 of hyperbolic type. This extension behaves near the origin as
 $r^{\beta + \alpha}$ and towards infinity as a power
 $r^{\beta}$. Since $r^{\beta + \alpha}$ is not infinitely often
 differentiable at zero, these functions provide a useful new class of
 approximants whose smoothness can be adjusted according to the
 smoothness of the function to be approximated. This makes them very
 flexible near zero because, as is the case with polynomial splines,
 we do not wish as a rule to approximate arbitrary approximands with
 an infinitely differentiable approximant as is the case with the
 usual radial basis functions -- hence their success.

As before, we write $g$ in a more convenient form
\begin{align*}g(x)&= r^{\beta} \left(\tanh r\right)^\alpha+r^{\beta}-r^{\beta}= r^{\beta} - r^{\beta} \left(1- \tanh^\alpha r \right)\\
&=  r^{\beta} - r^{\beta} \left(1- \left(\dfrac{\mathrm{e}^{r}-\mathrm{e}^{-r}}{\mathrm{e}^{r}+\mathrm{e}^{-r}} \right)^\alpha  \right)= 
 r^{\beta} - r^{\beta} \left(\dfrac{(\mathrm{e}^{r}+\mathrm{e}^{-r})^\alpha  - (\mathrm{e}^{r}-\mathrm{e}^{-r})^\alpha}{(\mathrm{e}^{r}+\mathrm{e}^{-r})^\alpha}  \right)\\
 &= r^{\beta} - r^{\beta}\dfrac{\mathrm{e}^{r\alpha}}{(\mathrm{e}^{r}+\mathrm{e}^{-r})^\alpha}\left((1+\mathrm{e}^{-2r})^\alpha - (1-\mathrm{e}^{-2r})^\alpha  \right)\\
  &= r^{\beta} - 2r^{\beta}\dfrac{\mathrm{e}^{r\alpha}}{(\mathrm{e}^{r}+\mathrm{e}^{-r})^\alpha}
 \sum_{k \; \rm odd}^{\infty}\dfrac{\Gamma(\alpha + 1)}{\Gamma(k + 1)\Gamma(\alpha - k + 1)}\mathrm{e}{-2rk}\\
   &= r^{\beta} - 2r^{\beta}\dfrac{1}{(1+\mathrm{e}^{-2r})^\alpha} \sum_{k \; odd}^{\infty}\dfrac{\Gamma(\alpha + 1)}{\Gamma(k + 1)\Gamma(\alpha - k + 1)}\mathrm{e}^{-2rk}\\
      &= r^{\beta} - 2r^{\beta}\left(\sum_{j=0}^{\infty}
 (-1)^j\mathrm{e}^{-2rj}\right)^\alpha \sum_{k \; \rm odd}^{\infty}\dfrac{\Gamma(\alpha + 1)}{\Gamma(k + 1)\Gamma(\alpha - k + 1)}\mathrm{e}^{-2rk}.
\end{align*}
We would now proceed similarly to the previous sections.

A parallel scheme can be followed if $g(x)= r^{\beta}\log r\left(\tanh r\right)^\alpha$ is considered instead.

\noindent {\it Remark.} It could also be considered $g(x)=r^{\beta} \tau (r)$, with $r^{\beta}(1-\tau (r))\in L^1(\mathbb{R}^n)$ or,  in a more general way,  $g(x)=\xi (r) \tau (r)$, with $\xi (r) (1-\tau (r))\in L^1(\mathbb{R}^n)$, and make the corresponding assumptions so that the Strang and Fix conditions are satisfied to achieve the associated polynomial reproduction degree and order of approximation by the quasi-interpolant. 

In fact, the more general case of a radial basis function $g$ with generalised Fourier transform
$$ \hat{g}(y)= G(s^2)+ \gamma s^{-\eta}, \quad y\in\mathbb{R}^n, \quad s:= \Vert y \Vert,$$
could be considered, by imposing adequate assumptions ($L^1$ integrability and other conditions) to construct the appropriate quasi-Lagrange function $\psi$ and the corresponding quasi-interpolant.
\section{Numerical experiments}
 Three graphical and numerical examples are now shown in order to illustrate the work. In the first example the trigonometric expansion associated with the quasi-Lagrange function transform is finite. In the second one, also in one dimension, a trigonometric series (but not a trigonometric polynomial) appears because the order of the singularity of the generalised Fourier transform of the chosen radial basis function is odd. The third example is two-dimensional and the expansion is, as in the first example, a trigonometric polynomial.
\subsection{Example 1}
Two quasi-interpolants will be compared: 
The first one is obtained from the radial basis function
$g_1(x)=(c^2+r^2)^p, \, x\in\mathbb{R}^n, \, c>0, \, p\in \mathbb{R}\setminus \mathbb{N}_0,  \, r:=\Vert x \Vert$, whose generalised Fourier transform is  (see \cite{gel}):
$$\hat{g_1}(y)= \frac{(2 \pi )^{\frac{n}{2}} 2^{p+1} c^{\frac{n}{2}+p} K_{\frac{n}{2}+p}(c \Vert y \Vert)}{\Gamma (-p)\Vert y \Vert^{\frac{n}{2}+p} }, \qquad y\in\mathbb{R}^n,$$
$K_\cdot$ being the modified Bessel function of the second kind, and $\Gamma$ the Gamma function.
The second quasi-interpolant comes from the radial basis function $g_2(x)=r^\beta \tanh(r)$, where $\beta\neq 2k$ and $\beta\neq -n-2k$ ($k=0,1,2,\ldots$). 
As we indicated in the previous sections, the generalised Fourier
transform needed to obtain the coefficients of the trigonometric
polynomial (and hence, the quasi-interpolant) is that of the function
$u(x)=r^\beta$, which is given by  (see the appendix) 
$$\hat{u}(y)= \frac{\pi ^{\frac{n}{2}} 2^{\beta +n} \Gamma
  \left(\frac{n+\beta }{2}\right) \Vert y\Vert^{-\beta -n}}{\Gamma
  \left(-\frac{\beta }{2}\right)},\qquad y\in\RR^n\setminus\{0\}.$$ 
We choose the values $n=1$,  $p=\frac32$ and $\beta= 3$. From \cite{olv} we read
\begin{align*}\label{Olver}
\frac{c^\nu}{\|y\|^{\nu}} K_\nu(c\|y\|) =& 2^{\nu-1}\frac{1}{\|y\|^{2\nu}}
\sum_{k=0}^{\nu-1} \dfrac{ (\nu-k-1)!}{k!(-4)^k}(c\|y\|)^{2k}+{}\\
&{}
-\left(\frac{-c^2}{2}\right)^\nu   \log \|y\|    \sum_{k=0}^{\infty}
\frac{(c\|y\|)^{2k}}{4^kk!\Gamma(\nu+ k+1)}  +{}
\\&{}-
\left(\frac{-c^2}{2}\right)^\nu\log c   \sum\limits_{k=0}^{\infty} \frac{(c\|y\|)^{2k}}{4^kk!\Gamma(\nu+ k+1)} +{}\\&+{} \left(\frac{-c^2}{2}\right)^\nu \sum_{k=0}^{\infty} \left(\frac{\log 2}{4^kk!\Gamma(\nu+
  k+1)}+\frac12\dfrac{\Psi(k+1)+\Psi(\nu+k+1)}{4^k(\nu+k)!k!}\right)\\
&{}\times(c\|y\|)^{2k},
\end{align*}
where $\nu\in\mathbb{N}$ and $\displaystyle\Psi(z)=\frac{\Gamma '(z)}{\Gamma (z)}$ is the Digamma function.  Therefore,
$$\hat{g_1}(y)= \frac{12}{y^4}-\frac{3 \pi ^2}{y^2}+ {\cal
  O}(1),\qquad y\to0.$$ 
On the other hand,
$$ \hat{u}(y)=  \frac{12}{y^4}.$$ 
This means that both generalised Fourier transforms present a singularity of order 4 at the origin.   We must cancel this singularity with the help of a suitable trigonometric polynomial so that the generalised Fourier transforms of the quasi-Lagrange functions possess the appropriate conditions to reproduce algebraic polynomials of degree 3, which is the maximum possible degree of reproduction.

We choose a $C^3(\mathbb{R})$ target function which is not $C^4(\mathbb{R})$, so that the approximand and the approximant of hyperbolic type have the same smoothness, which is desirable as mentioned above. 

Thus, we take 
\[ f(x)=\left(1-x^2\right)_+^4.
\]

Now we choose the ordered set of points $$\Delta= \left\{j-5 :\, j= 1, \ldots, 9\right\}.$$ 

Our goal is to solve a linear system of equations $Ax=b$ where $A=(a_{i,j})_{1\leq i,j\leq 9}$ being $a_{i,j}=(j-5)^{i-1}.$ The nine components of the solution will be the coefficients of the trigonometric polynomial $P(y)= \displaystyle\sum_{k=-4}^4 \mu_k \mathrm{e}^{-ik y}$ that we are looking for.

Although the matrix is the same for both cases, the expression of the independent term of the system varies depending on whether we are working with $g_1$ or $u$. In both cases, the solution of the system must satisfy 
$$P({y}) \hat g_1({y})=1+{\cal O}(\Vert {y} \Vert ^{4}) \quad \mbox{and} \quad P({y}) \hat u ({y})=1+{\cal O}(\Vert {y} \Vert ^{4}).$$
In particular, if we work with $\hat{g_1}(y)$ the
independent term $b\in\mathbb{R}^9$  will have all its components null except 
$$b(5)=2, \quad b(7)=-15c^2, \quad b(9)=\frac{105}{2} c^4 (4 \log (c)+4 \gamma +1-4 \log (2)),$$ whereas 
if we work with $\hat u (y)$ the only non-zero component will be $b(5)= 2. $

Taking $c=0. 5$,  the values $\displaystyle \{\mu_{1,k}\}_{k=-4}^4$,
solution vector of the linear system for the first case, can be found
in \cite{jjmm},  whereas for the second one are
$$\displaystyle \{\mu_{2,k}\}_{k=-4}^4=\left\{\frac{7}{2880}, -\frac{1}{30},\frac{169}{720},-\frac{61}{90},\frac{91}{96},-\frac{61}{90},\frac{169}{720},-\frac{1}{30},\frac{7}{2880}\right\}. $$
Now, we can construct both quasi-Lagrange functions, whose expression is
\begin{equation*}
 \psi_i(x)=\sum_{k\in\Delta} \mu_{i,k} \, g_i(x-k),\quad x\in\mathbb{R}, \quad i=1,2.
\end{equation*}
and its associated quasi-interpolants
\begin{equation*}
Q_if(x)=\sum_{j\in\mathbb{Z}}f(hj)\psi_i(x/h-j),\quad x\in\mathbb{R}, \quad i=1,2.
\end{equation*}
We now show some graphics and the maximum error obtained. We have taken $h=10^{-3}.$
\begin{figure}[h]\label{fig1}  
    \centering
    \begin{minipage}{0.45\textwidth}
        \centering
        \includegraphics[width=0.9\textwidth]{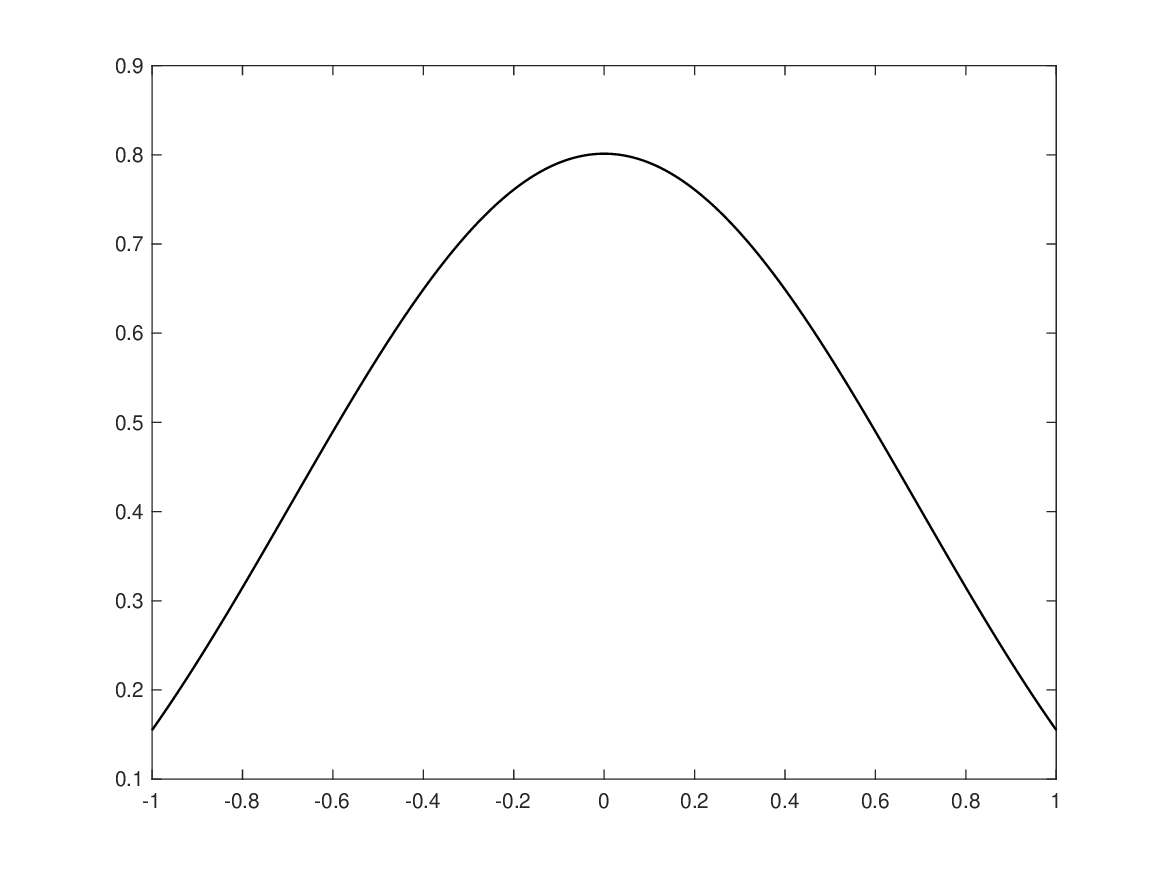} 
    \end{minipage}\hfill
    \begin{minipage}{0.45\textwidth}
        \centering
        \includegraphics[width=0.9\textwidth]{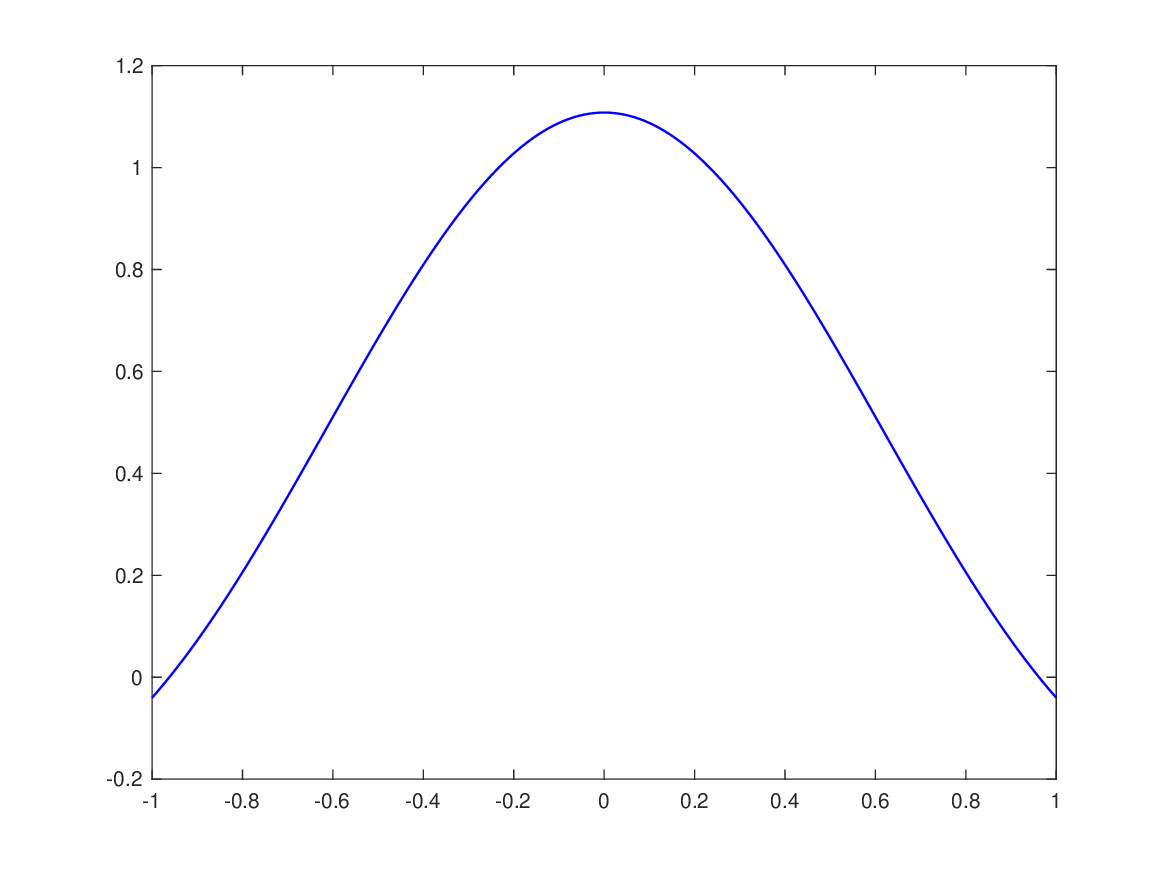}
    \end{minipage}
 \caption{Quasi-Lagrange functions corresponding to $g_1$ (left) and  $g_2$ (right).}
\end{figure}

 \begin{figure}[h]\label{fig2}  
    \centering
    \begin{minipage}{0.45\textwidth}
        \centering
        \includegraphics[width=0.9\textwidth]{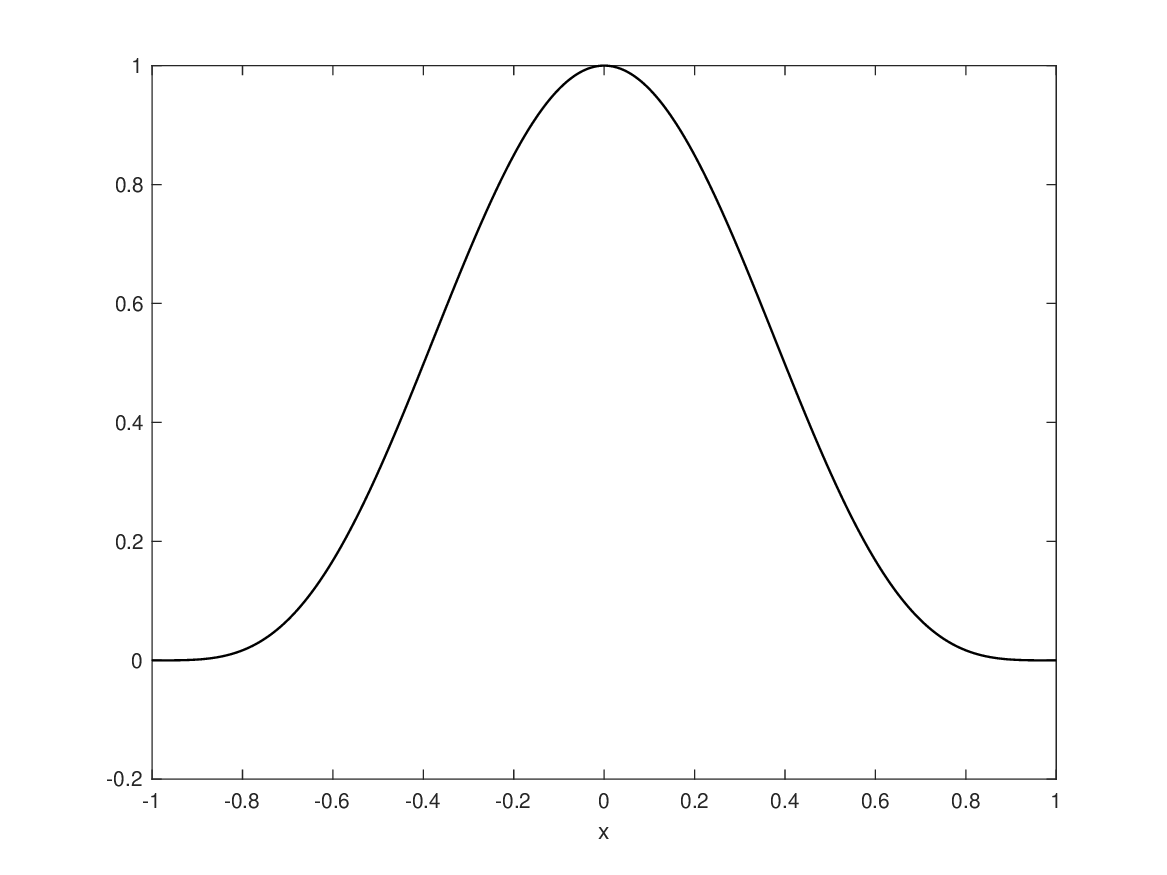} 
    \end{minipage}\hfill
    \begin{minipage}{0.45\textwidth}
        \centering
        \includegraphics[width=0.9\textwidth]{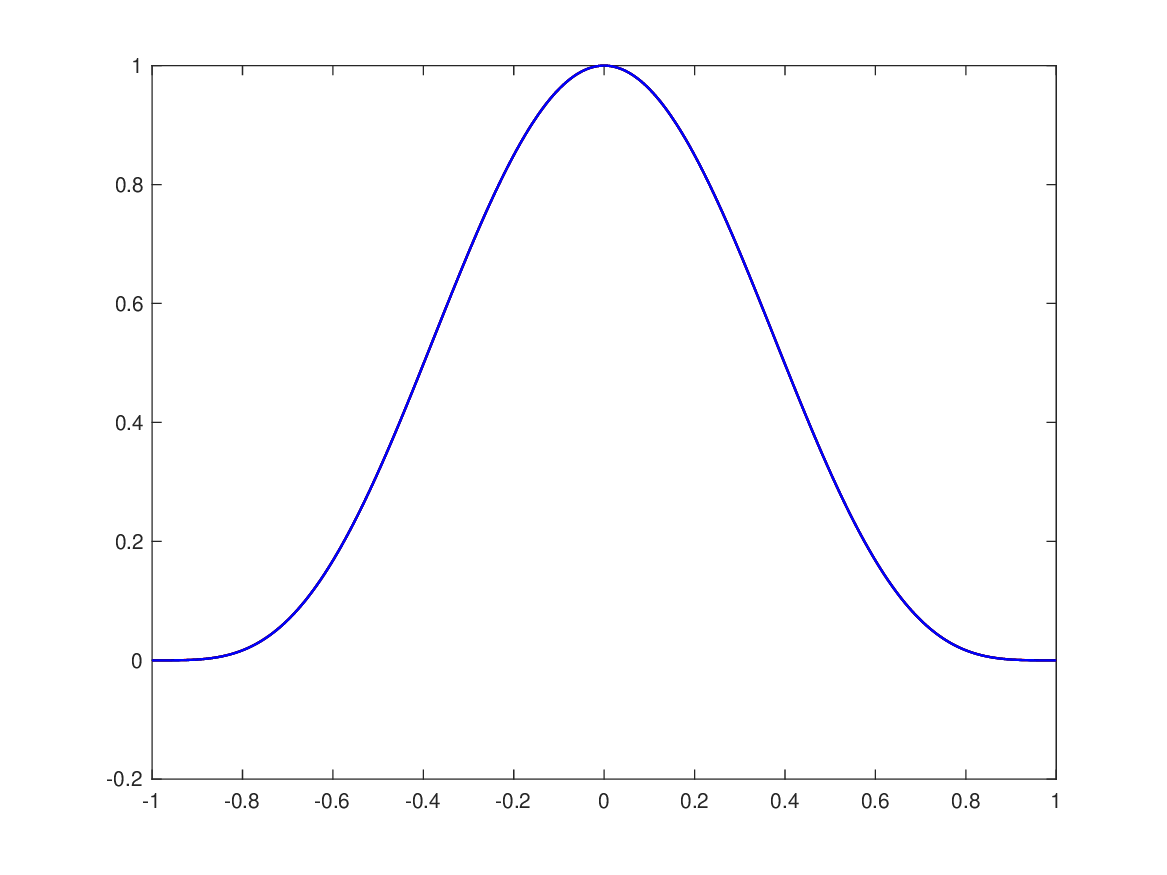}
    \end{minipage}
 \caption{Quasi-interpolants corresponding to $g_1$ (left) and  $g_2$ (right).}
\end{figure}

 \begin{figure}[!h]\label{fig3}  
        \centering
    \begin{minipage}{0.45\textwidth}
        \centering
        \includegraphics[width=0.9\textwidth]{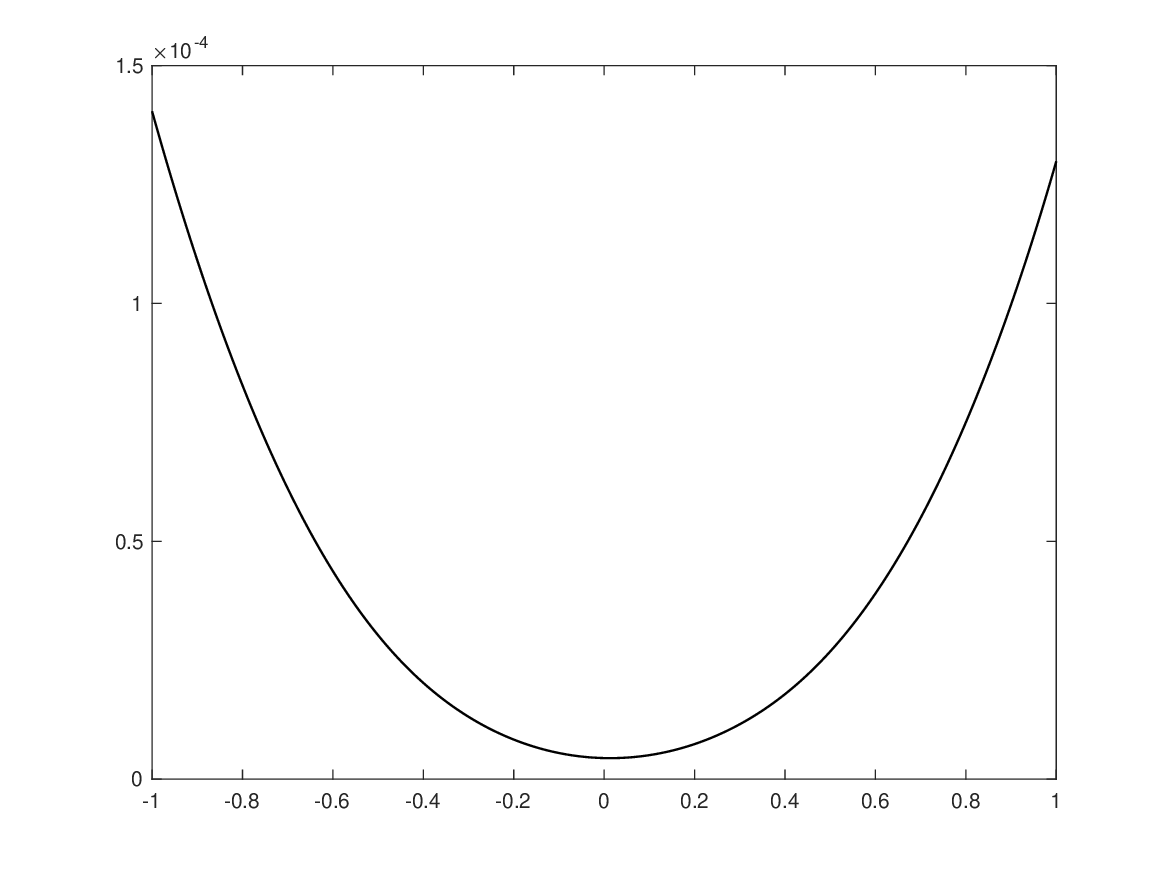} 
    \end{minipage}\hfill
    \begin{minipage}{0.45\textwidth}
        \centering
        \includegraphics[width=0.9\textwidth]{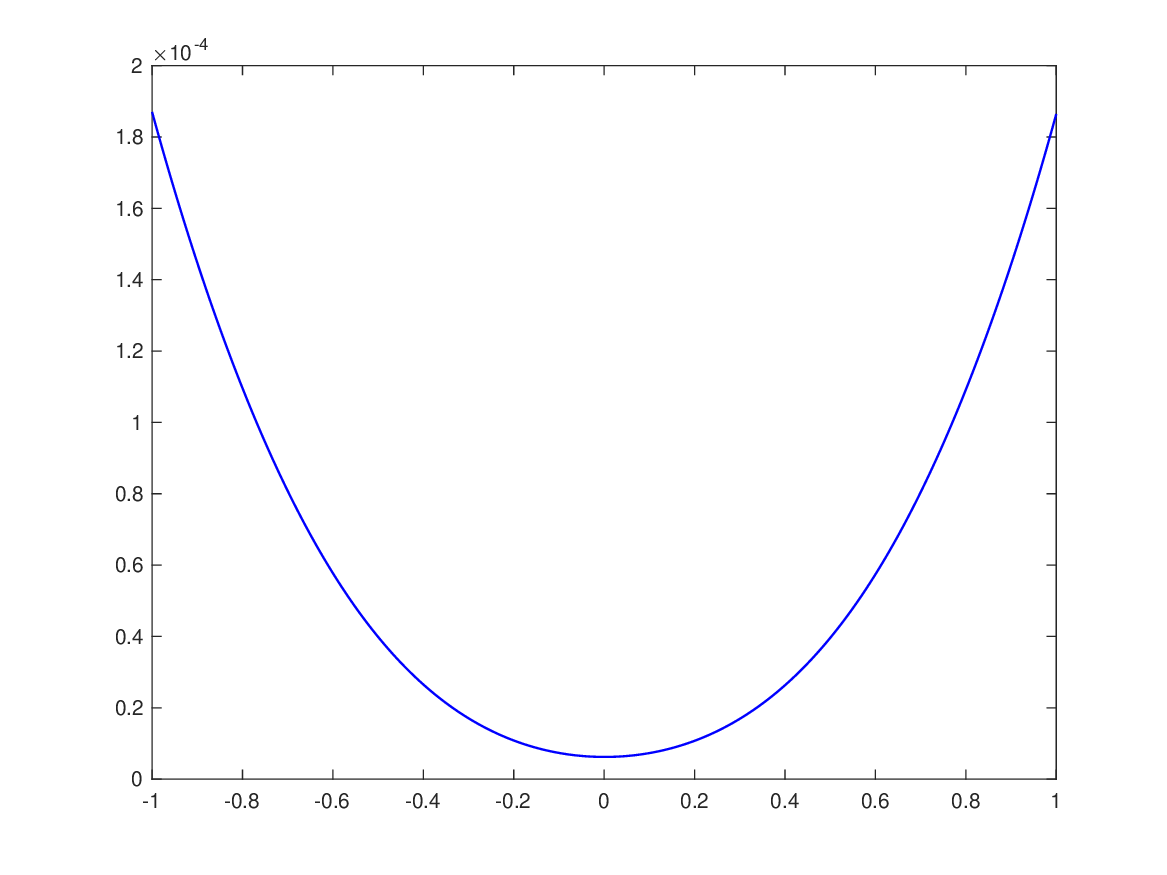}
    \end{minipage}
 \caption{Graphics of the  absolute error between  $f(x)$ and the quasi-interpolant based on  $g_1(x)$ (left)  and the one based on  $g_2(x)$  (right).}
 \end{figure}

\begin{table}[!h]
\centering
\begin{tabular}{|c|c|c|c|}
\hline
RBF   & Sing. at 0 & Polyn. reproduction &  Max. error  \\ \hline
 $g_1(x)=(c^2+r^2)^{3/2}$   &-4 &3  & 0.0001404 \\ \hline
 $g_2(x)=r^3 \tanh(r)$    &-4& 3 &  0.000187 \\ \hline
\end{tabular}
\caption{Summary of Example 1}
\label{tab1}
\end{table}
Although the smoothness of the target function and  $g_2$ are the same (they are $C^3$ functions), we obtain a slightly better approximation with the quasi-interpolant based on the generalized multiquadric, $g_1$.
\subsection{Example  2}
The radial basis functions for this second example are
$$g_1(x) =(c^2+r^2) \log \left(\sqrt{c^2+r^2}\right)-c^2 \log (c)-r^2 \log (c)  $$ and
$$g_2(x) =\left(r^2\log(r)+ \gamma r^2 \right) \tanh (r),$$
where $x\in\mathbb{R}^n$, $c>0$, $r:=\Vert x \Vert$,
and $\gamma$ is the Euler-Mascheroni constant.  The generalised Fourier transform of $g_1$ can be found in \cite{ort3} and is given by 
$$\hat {g_1}(y)=2^{\frac{n}{2}+1} \pi ^{\frac{n}{2}} \left(\frac{c}{\Vert y\Vert}\right)^{\frac{n}{2}+1} K_{\frac{n}{2}+1}(c \Vert y\Vert).$$
Let us consider the one-dimensional case, that is, $n=1$. Therefore, 
$$\hat {g_1}(y)=2^{\frac{3}{2}} \pi ^{\frac12} \left(\frac{c}{\vert y\vert}\right)^{\frac{3}{2}} K_{\frac{3}{2}}(c \vert y\vert).$$
Taking into account that (see \cite[10.2.17]{abr})
$$K_{\frac{3}{2}}(c \vert y\vert)=\dfrac{\mathrm{e}^{-c \vert y\vert} (\pi/2)^{1/2} \left(1+\frac{1}{c \vert y\vert}\right)}{(c \vert y\vert)^{1/2}},$$
we have that
$$\hat{g_1}(y)=2\pi \mathrm{e}^{-c \vert y \vert} \left( \frac{1}{\vert y \vert ^3}+\frac{c}{\vert y \vert^2}\right),$$
whose expansion about the origin is given by
$$2\pi\left( \frac{1}{\vert y\vert^3}-\frac{c^2}{2\vert y \vert} +\frac{c^3}{3}-\frac{c^4}{8}\vert y \vert+\cdots\right).$$
The generalised Fourier transform of $u(x)=r^2\log(r)+ \gamma r^2$ is
given by (see the appendix)
$$\hat{u}(y)=\dfrac{2\pi}{\vert y\vert^3}.$$
In both cases, there exists a singularity of order 3 at the
origin. Therefore, it will be needed to have a trigonometric series (but not a
trigonometric polynomial) to annihilate this singularity.  The
coefficients of this trigonometric expansion will be computed as we
did in  \cite{jjmm}, i.e., these coefficients will be, for instance,
the Fourier coefficients of the function $\dfrac{1}{2\pi}(2-2
\cos(x))^{3/2}$. As there are infinitely many coefficients, we need to
truncate them in order to obtain the graphics and the error we
present. For this case we have considered $c=0.5$ and $2^{11}$ Fourier
coefficients (the convergence of the corresponding Fourier series is
very slow) and for the quasi-interpolants we have chosen  $h=10^{-2}.$
The target function will be  
\[ f(x)=\left(1-x^2\right)_+^3.
\]
We observe that $f\in C^2(\mathbb{R})\setminus C^3(\mathbb{R}), $ analogously to  $g_2$. 

 \begin{figure}[!h]\label{fig1ej2}  
    \centering
    \begin{minipage}{0.45\textwidth}
        \centering
        \includegraphics[width=0.9\textwidth]{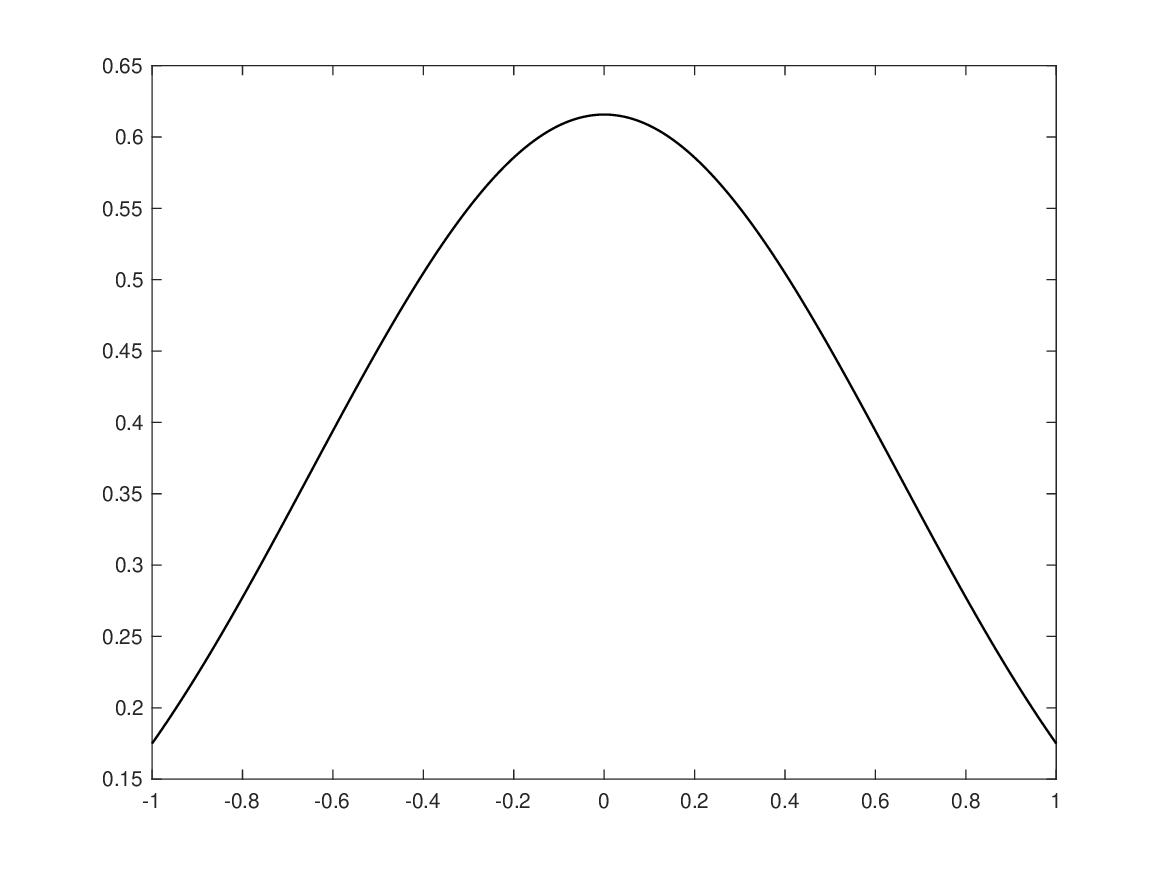} 
    \end{minipage}\hfill
    \begin{minipage}{0.45\textwidth}
        \centering
        \includegraphics[width=0.9\textwidth]{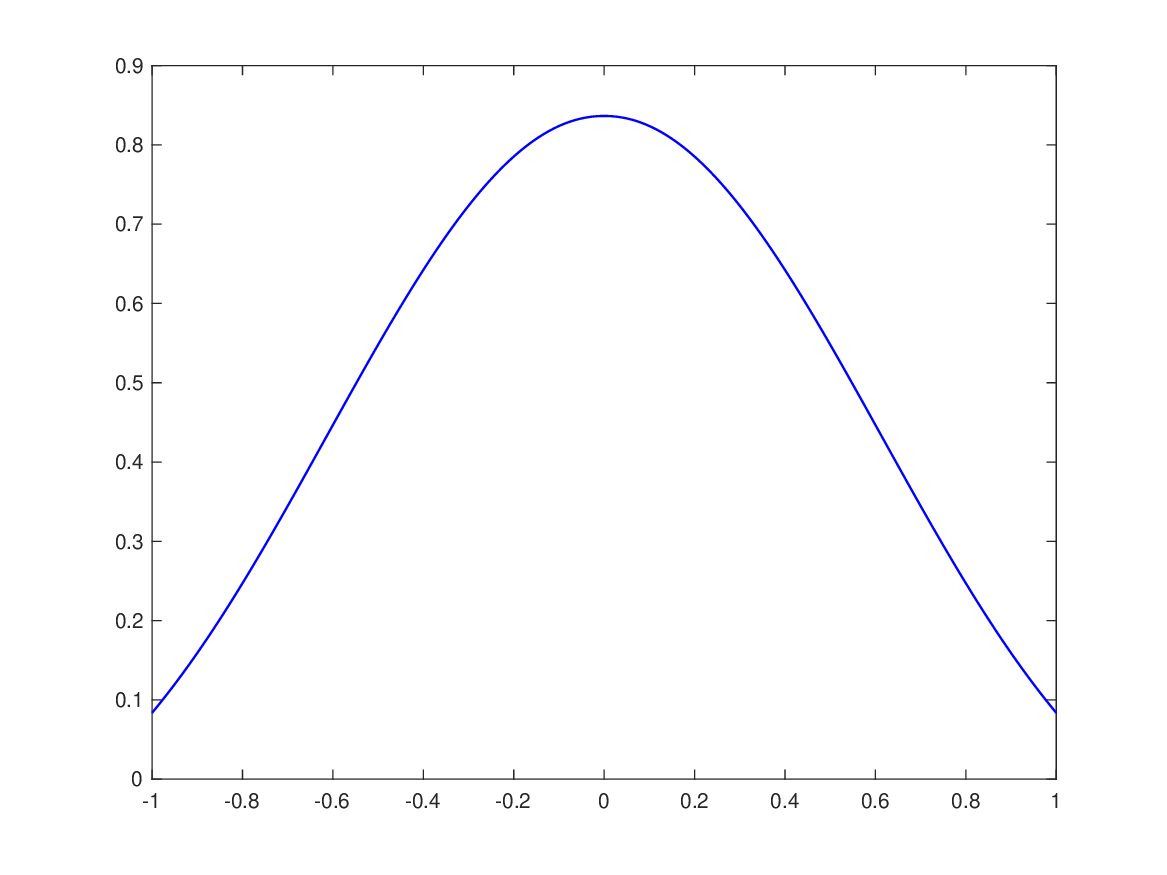}
    \end{minipage}
 \caption{Quasi-Lagrange functions corresponding to $g_1$ (left) and  $g_2$ (right).}
\end{figure}

 \begin{figure}[!h]\label{fig2ej2}  
    \centering
    \begin{minipage}{0.45\textwidth}
        \centering
        \includegraphics[width=0.9\textwidth]{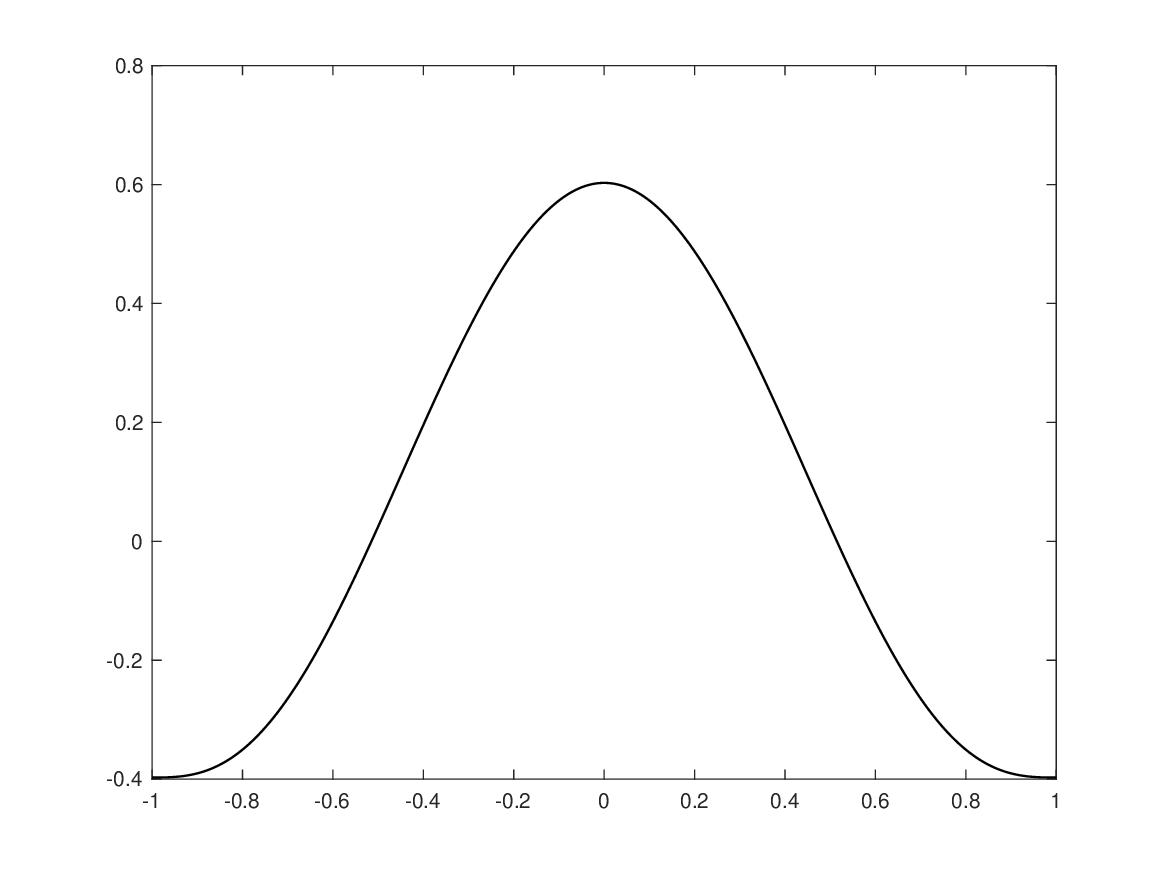} 
    \end{minipage}\hfill
    \begin{minipage}{0.45\textwidth}
        \centering
        \includegraphics[width=0.9\textwidth]{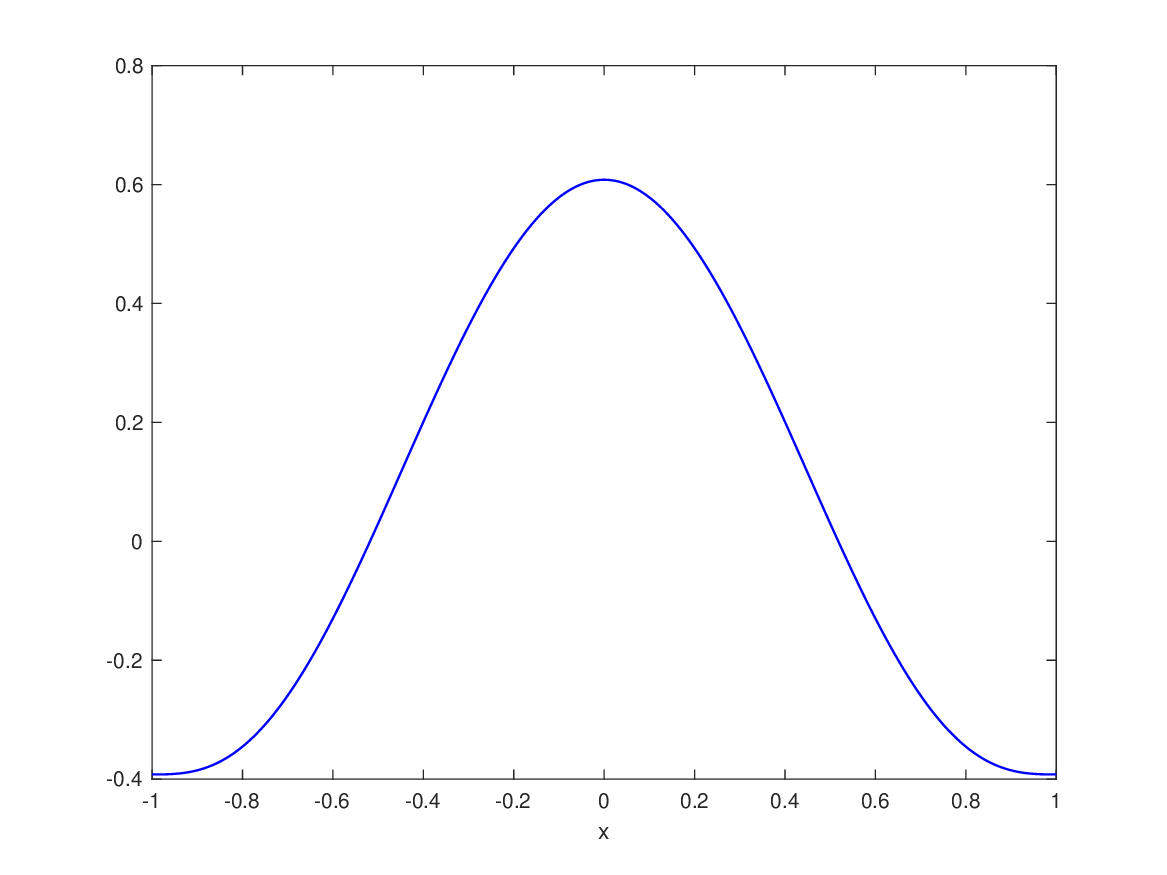}
    \end{minipage}
 \caption{Quasi-interpolants corresponding to $g_1$ (left) and  $g_2$ (right).}
\end{figure}

 \begin{figure}[!h]
        \centering
    \begin{minipage}{0.45\textwidth}
        \centering
        \includegraphics[width=0.9\textwidth]{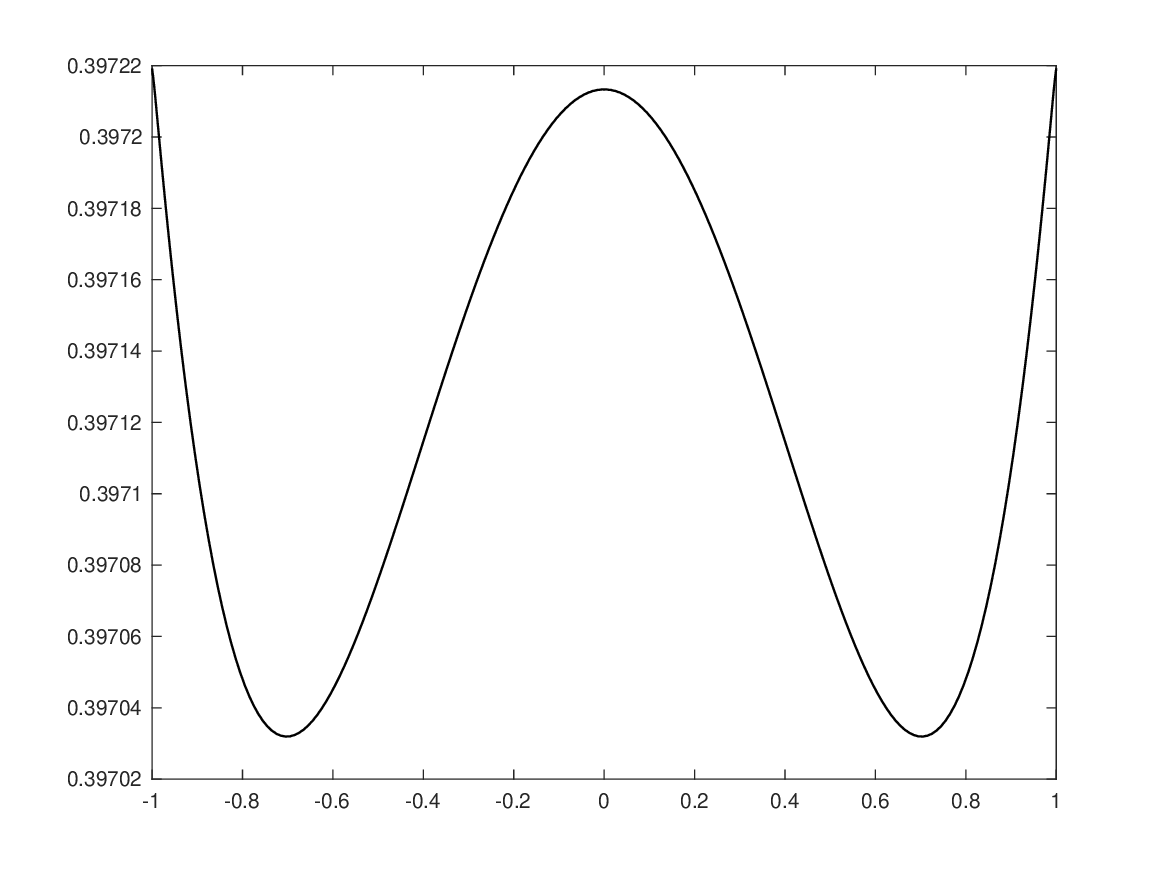} 
    \end{minipage}\hfill
    \begin{minipage}{0.45\textwidth}
        \centering
        \includegraphics[width=0.9\textwidth]{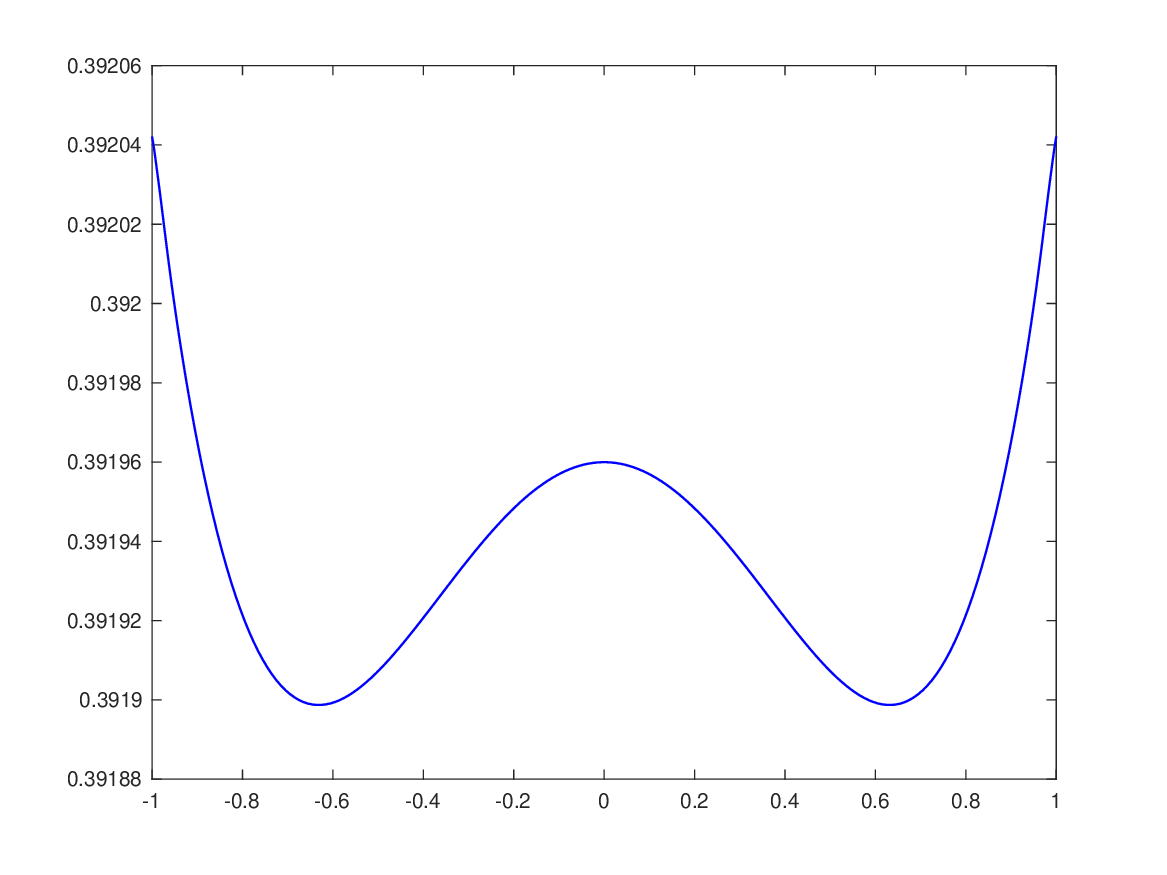}
    \end{minipage}
 \caption{Graphics of the  absolute error between  $f(x)$ and the quasi-interpolant based on  $g_1(x)$ (left)  and the one based on  $g_2(x)$  (right).}
 \label{fig3ej2}  
 \end{figure}

\begin{table}[!h]
\centering
\begin{tabular}{|c|c|c|c|}
\hline
RBF   & Sing. at 0 & Polyn. reproduction &  Max. error  \\ \hline
 $g_1  $   &-3 &1  & 0.39722 \\ \hline
 $g_2 $   &-3& 1 &   0.39204 \\ \hline
\end{tabular}
\caption{Summary of Example 2}
\label{tab2}
\end{table}
By looking at  Table \ref{tab2} and Figure \ref{fig3ej2}  we note again that the approaches are very similar.
\subsection{Example 3}
This final example is developed in two dimensions, that is, $n=2$.
We compare the quasi-interpolants obtained from
$$g_1(x) =(c^2+r^2) \log (\sqrt{c^2+r^2})-c^2\log (c)-r^2 \log (c)  $$
and
$$g_2(x) =\left(r^2\log(r)+ (\gamma- \log(2)) r^2 \right) \tanh (r).$$ 
The generalised Fourier transform of  $g_1$ is
$$\hat{g_1}(y)=\dfrac{4 \pi c^2}{\Vert y\Vert^2} K_{2}(c \Vert y\Vert)$$
and the one of 
 $u(x)=r^2\log(r)+ (\gamma- \log(2)) r^2$ is (see the appendix)
 $$\hat{u}(y)= \dfrac{8 \pi}{\Vert y\Vert^4}.$$
Both of them have a singularity of order -4, which can be annihilated by a suitable trigonometrical polynomial.
For this task, we will need a set of 21 points near the origin, in order to impose the Strang-Fix conditions. We set 
\begin{align*}
\Delta=&\left\{(0,0),(1,0),  (0,1),(-1,0),(0,-1), (1,1),(-1,1), (-1,-1), (1,-1), (2,0), (0,2), \right\} \\
  &\left\{ (-2,0), (0,-2), (3,0),(0,3), (0,-3), (-3,0), (2,2), (-2,2), (-2,-2), (2,-2)\right\}
\end{align*}
These number of points, 21, is sufficient for obtaining a solution of the linear system $Ax=b.$  This solution will be the coefficients for generating  both quasi-Lagrange functions.
In this case, $A \in \mathcal{M}_{36\times21}.$ Each row is developed as follows: we take the canonical basis of $\mathbb{P}_7$ in two variables: $\{1, x, y, x^2, xy, y^2, \ldots, x^2y^5, xy^6, y^7\}.$  Its dimension is $\binom{2+7}{2}=36.$ Now, each row of $A$ is obtained by evaluating each point of $\Delta$ on the basis element, i.e, the first row is composed by ones, the second one by the first component of the points, etc.

When working with $g_1$ the independent term, $b,$  has null components  except $$b(11) = b(15)=\frac{3}{\pi}, \quad b(13) =\frac{1}{\pi}, \quad b(22) = b(28)=-\frac{45c^2}{2},\quad  b(24) = b(26)= -\frac{9c^2}{2\pi},$$
and for $u$ the unique nonzero components of $b$ are
$$b(11)=\frac{3}{\pi}, \quad b(13) =\frac{1}{\pi}, \quad  b(15)=\frac{3}{\pi}.$$ 
In order two obtain the quasi-Lagrange functions and the quasi-interpolants we have set $c=0.5$ and $h=0.01.$ 
  \begin{figure}[!h]\label{fig0ej3}  
        \centering
        \includegraphics[width=0.5\textwidth]{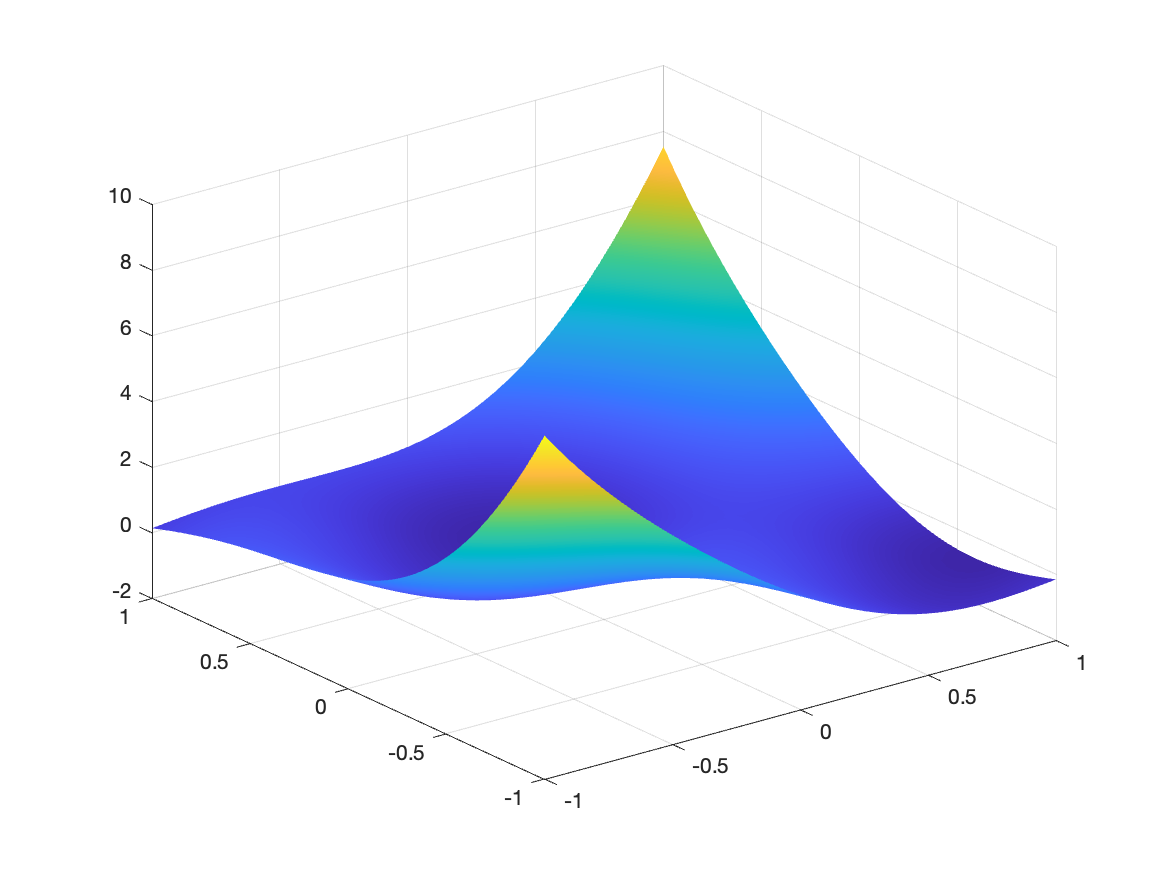} 
 \caption{Target function.}
 \end{figure}
 
  \begin{figure}[!h]\label{fig1ej3}  
    \centering
    \begin{minipage}{0.45\textwidth}
        \centering
        \includegraphics[width=\textwidth]{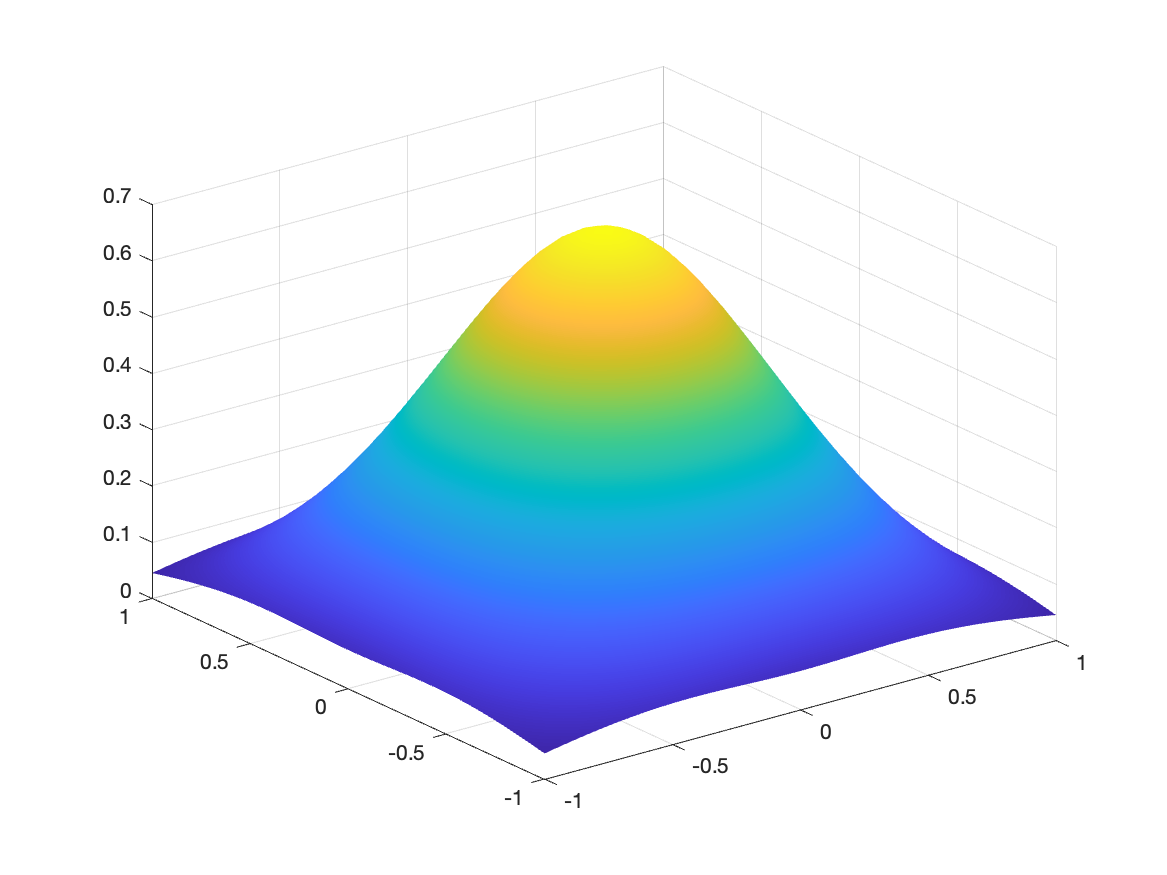} 
    \end{minipage}\hfill
    \begin{minipage}{0.45\textwidth}
        \centering
        \includegraphics[width=\textwidth]{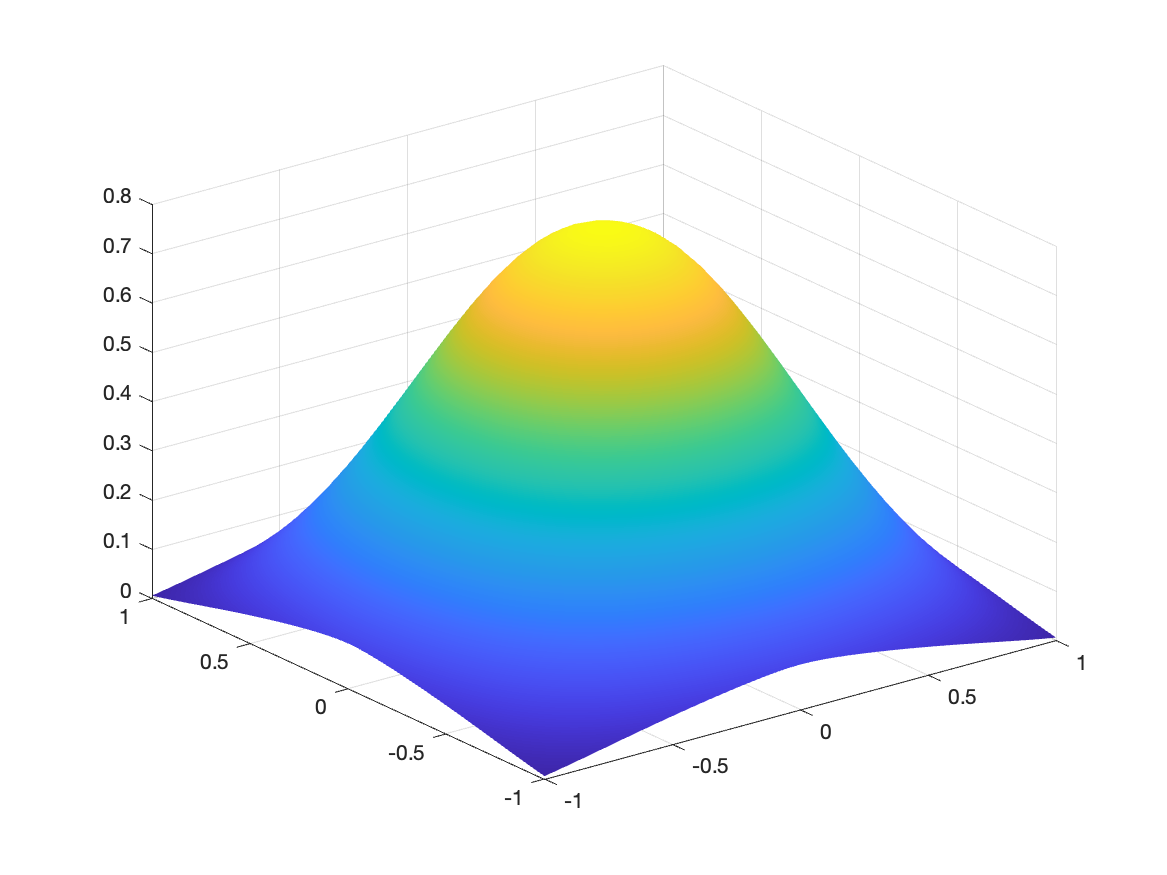}
    \end{minipage}
 \caption{Quasi-Lagrange functions corresponding to $g_1$ (left) and  $g_2$ (right) with $c=0.5$.}
\end{figure}
The chosen target function is $f(x_1,x_2)=(x_1 + x_2)^2 |x_1 + x_2| + \cos(2 x_1 - x_2) \sin(x_1 - 2 x_2)\in C^2(\mathbb{R}^2)\setminus C^3(\mathbb{R}^2),$   analogously to $g_2$ again.

\begin{figure}[h]\label{fig2ej3}  
    \centering
    \begin{minipage}{0.45\textwidth}
        \centering
        \includegraphics[width=0.9\textwidth]{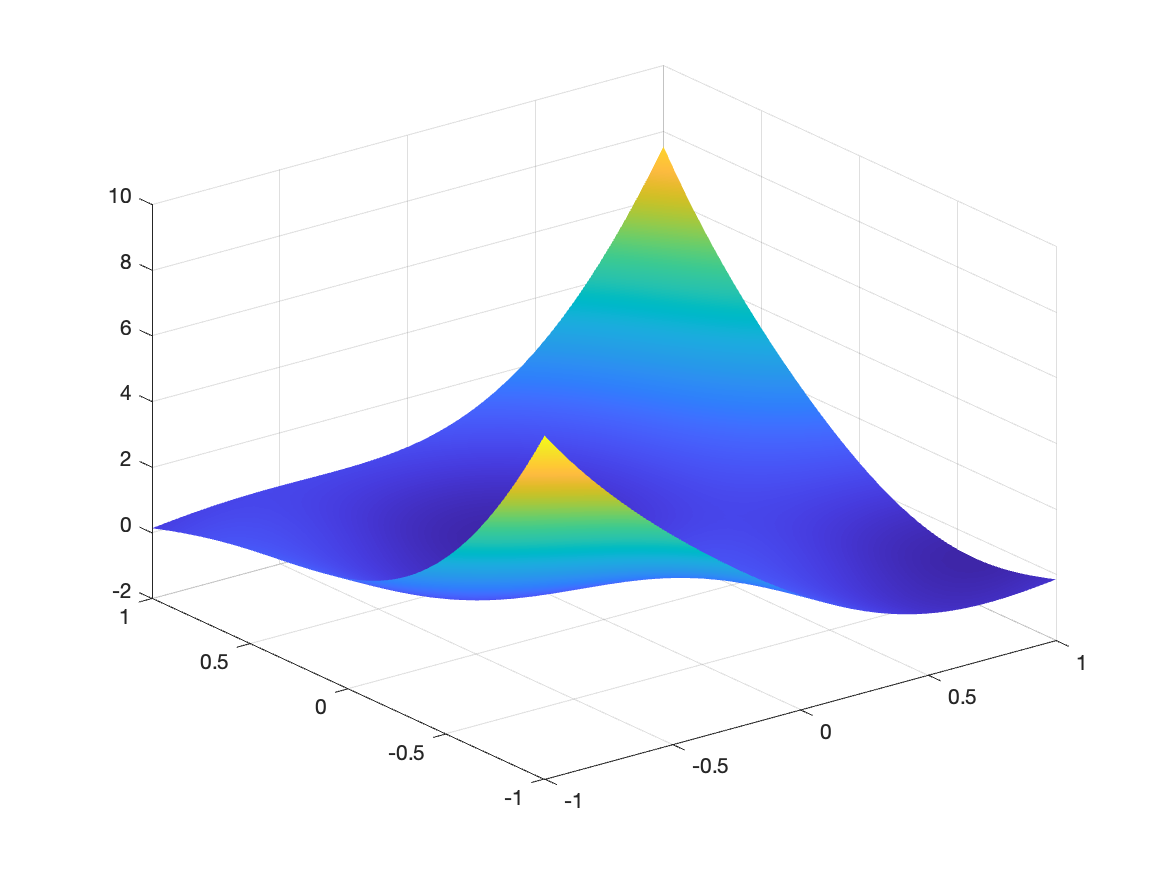} 
    \end{minipage}\hfill
    \begin{minipage}{0.45\textwidth}
        \centering
        \includegraphics[width=0.9\textwidth]{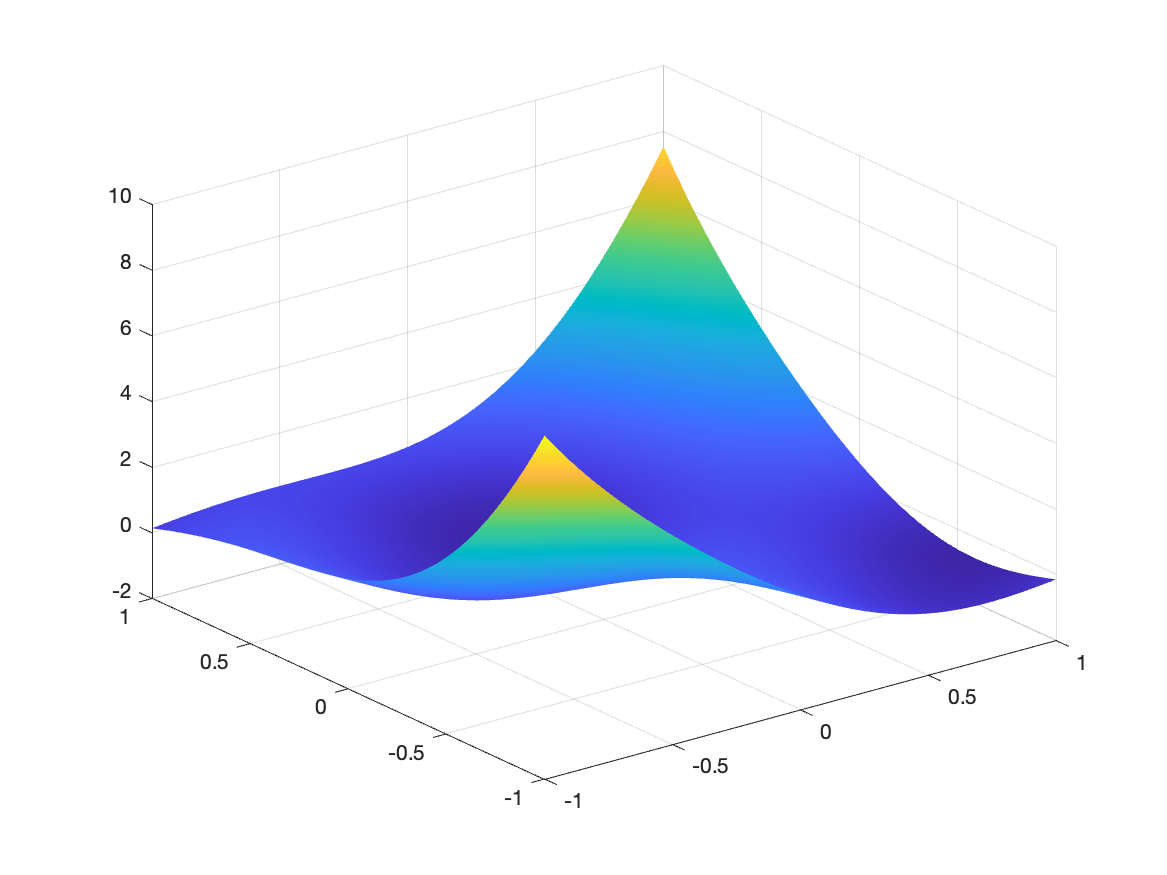}
    \end{minipage}
 \caption{Quasi-interpolants corresponding to $g_1$ (left) and  $g_2$ (right).}
\end{figure}

\begin{figure}[h]\label{fig3ej3}  
    \centering
    \begin{minipage}{0.45\textwidth}
        \centering
        \includegraphics[width=0.9\textwidth]{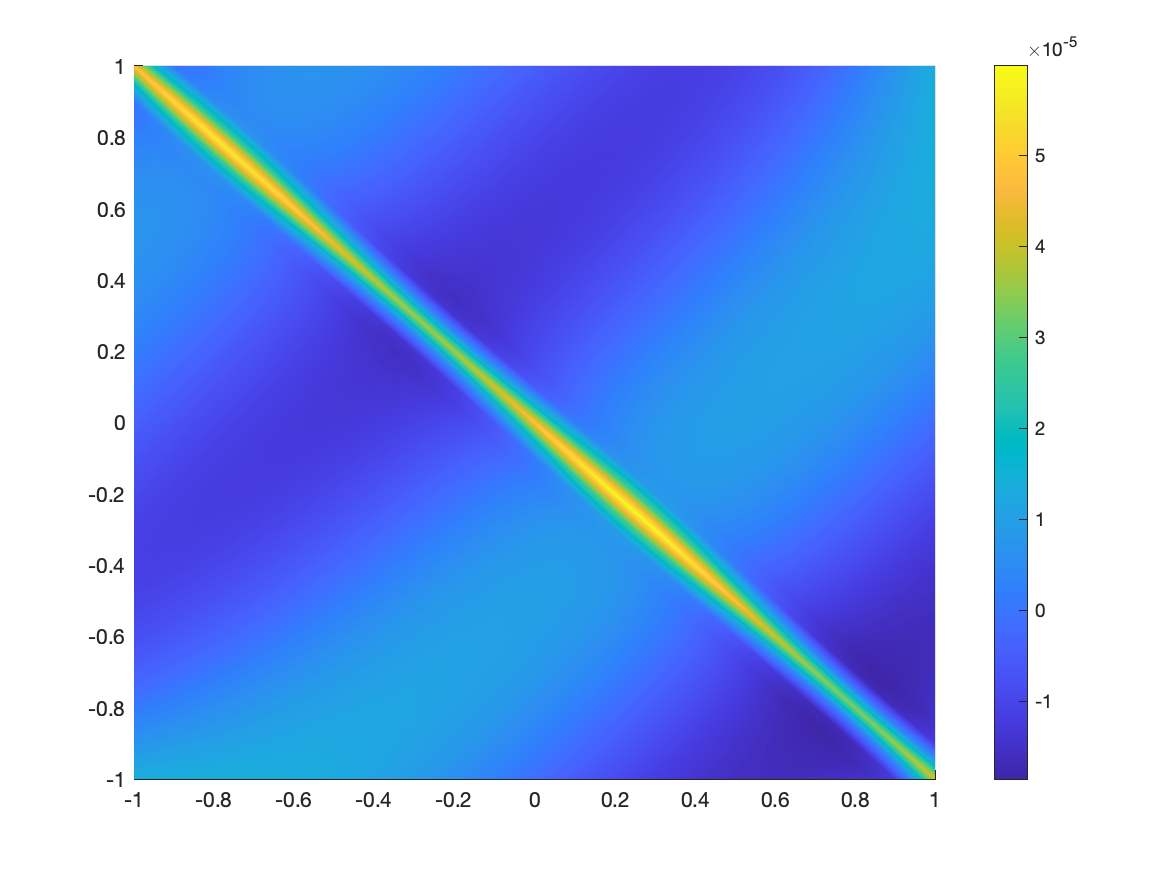} 
    \end{minipage}\hfill
    \begin{minipage}{0.45\textwidth}
        \centering
        \includegraphics[width=0.9\textwidth]{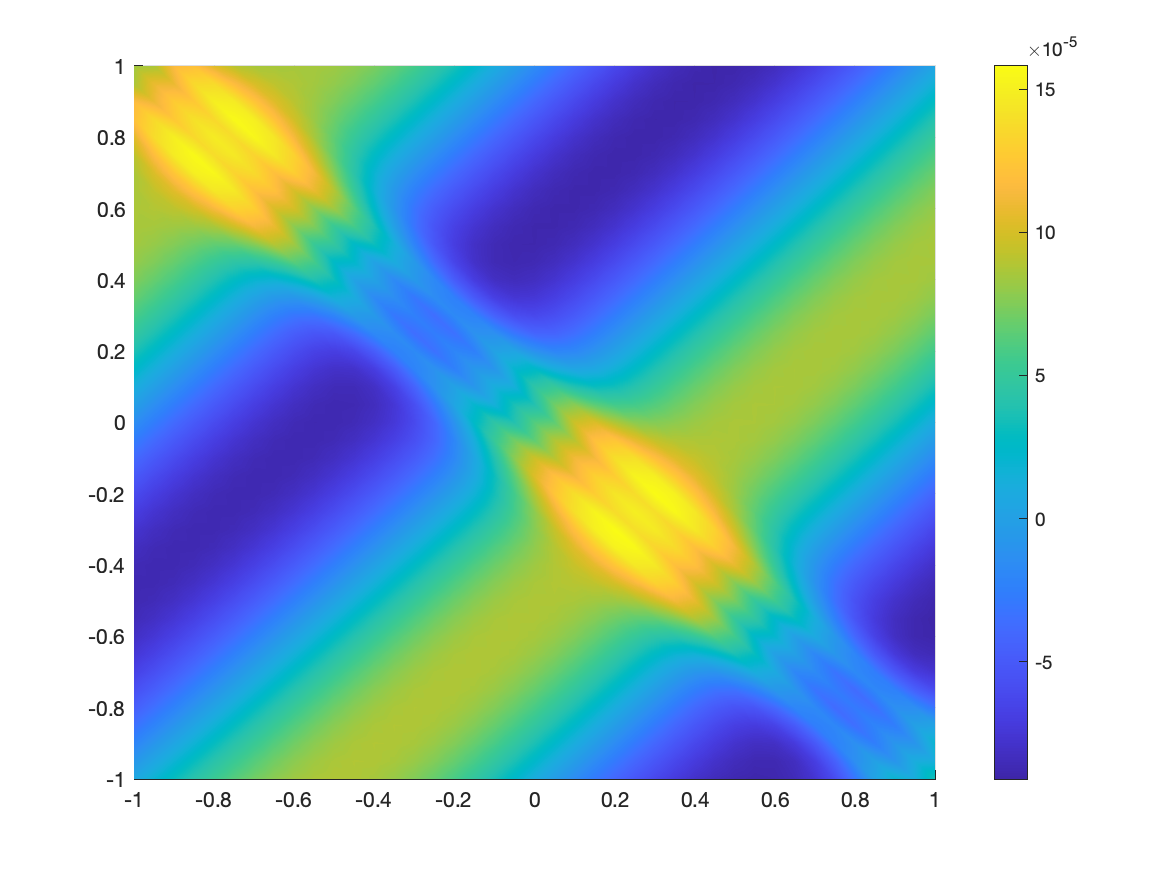}
    \end{minipage}
 \caption{Graphics of the  absolute error between  $f(x)$ and the quasi-interpolant based on  $g_1(x)$ (left)  and the one based on  $g_2(x)$  (right).}
\end{figure}

 \begin{table}[!h]
\centering
\begin{tabular}{|c|c|c|c|c|}
\hline
RBF   & Sing. at 0 & Polyn. reproduction &  Max. error & RMSE  \\ \hline
 $g_1  $   &-4 &3  &     0.000060 &0.000012\\ \hline
 $g_2 $   &-4& 3 &  0.000068 &  0.000158  \\ \hline
\end{tabular}
\caption{Summary of Example 3}
\label{tab3}
\end{table}

\section*{Acknowledgements}
This work has been partially financed by Junta de Andaluc{\'\i}a (Research groups FQM178 and FQM191) and University of Ja\'en (Research group EI\_FQM8). 

\section*{Appendix}
For the sake of completeness, the explicit univariate generalised Fourier transform of the hyperbolic tangent and the generalised Fourier transforms of some classical radial basis functions are included in this section.

\subsection*{A.1. Univariate generalised Fourier transform of $\text{tanh}$}

Taking into account that $$\frac{\rm{d}}{\rm{d}t}\text{tanh} (t)=\text{sech}^2 (t),$$ 
we have that
 $$\rm{i} \omega \mathcal{F} \left(\text{tanh} \right)(\omega)=\mathcal{F}   \text{sech}^2 (\omega).$$ 

Now, let us  compute 
$$\mathcal{F}   \text{sech}^2 (\omega)=\int_{-\infty}^{\infty} \text{sech}^2 (t ){\rm e}^{-\rm{i} t \omega}\rm{d}t.$$
This will be done in two steps:

{\bf Step 1.} We will  compute $\mathcal{F}    \text{sech} (\omega)$.
First, for  $0 < \operatorname{Re} \alpha < \operatorname{Re} \beta$, by the expression 3.241 (2) of (\cite{grad}), we have  that
\begin{equation*} \int \limits_0^\infty \frac{x^{\alpha -1}}{1+x^\beta} \, {\rm d} x = \frac{\pi}{\beta} \csc\left(\frac{\alpha}{\beta} \pi \right) \, .
\end{equation*}
Now, we can compute
\begin{align*}
\mathcal{F}   \text{sech} (\omega)&=\int\limits_{-\infty}^\infty
\frac{\mathrm{e}^{-\mathrm{i} \omega x}}{\cosh(x)} \, \mathrm{d} x =
\left\{t=\mathrm{e}^{x}\right\}= 2 \int \limits_0^\infty \frac{t^{-\mathrm{i}
    \omega}}{1+t^2} \, \mathrm{d} t = 2 \frac{\pi}{2} \csc
\left(\frac{1 - \mathrm{i} \omega}{2} \pi \right) =  
\pi \sec \left(\mathrm{i} \frac{\pi\omega}{2} \right) \\ &= 
\pi \text{sech} \left(- \frac{\pi \omega}{2}\right) \, =
 \pi \text{sech} \left(\frac{\pi \omega}{2}\right) \, .
\end{align*}

{\bf Step 2.}
We compute $\mathcal{F}   \text{sech}^2 (\omega).$ We use the property
\begin{equation}\label{tanh}
\mathcal{F}(f^2)=\frac{1}{(2\pi) ^n}\mathcal{F}(f)*\mathcal{F}(f),
\end{equation}
and taking $f(t)  =\text{sech} (t)$   the right hand side of (\ref{tanh}) is
\begin{align*}
I(\omega)&= \frac{1}{2\pi}\int_{-\infty}^{\infty}\pi \text{sech}\,\left(\frac{\pi t}{2}\right) \pi \text{sech}\,\left(\frac{\pi(\omega-t)}{2}\right)\,{\rm d} t= \frac{\pi}{2}\int_{-\infty}^{\infty}  \text{sech}\,\left(\frac{\pi t}{2}\right)\text{sech}\,\left(\frac{\pi(\omega-t)}{2}\right) {\rm d} t\\ 
&=2\pi \int_{-\infty}^{\infty}  \dfrac{\mathrm{e}^{-\frac{\pi t}{2}}}{1+\mathrm{e}^{-\pi t}} \dfrac{\mathrm{e}^{-\frac{\pi(\omega-t)}{2}}}{1+\mathrm{e}^{-\pi(\omega-t)}}\,{\rm d} t 
= 2\pi \mathrm{e}^{-\frac{\pi \omega}{2}} \int_{-\infty}^{\infty}  \dfrac{1}{\left(1+\mathrm{e}^{-\pi t}\right)\left(1+\mathrm{e}^{-\pi \omega}\mathrm{e}^{\pi t} \right)} \,{\rm d} t
\\ &=2\pi \mathrm{e}^{\frac{\pi \omega}{2}}  \int_{-\infty}^{\infty}  \dfrac{ 1}{\left(1+\mathrm{e}^{-\pi t}\right)\left(\mathrm{e}^{\pi \omega}+\mathrm{e}^{\pi t} \right)}\,{\rm d} t
=  -2 \mathrm{e}^{\frac{\pi \omega}{2}} \int_{\infty}^{0}\dfrac{1}{(1+u)(\mathrm{e}^{\pi \omega}+ \frac{1}{u})} \dfrac{1}{u}{\rm d} u \\
&= 2 \mathrm{e}^{\frac{-\pi \omega}{2}} \int^{\infty}_{0}\dfrac{1}{(1+u)(\mathrm{e}^{-\pi \omega}+ u)}{\rm d} u
= \dfrac{2\mathrm{e}^{\frac{-\pi \omega}{2}}}{\mathrm{e}^{-\pi \omega}-1} \int^{\infty}_{0} \left(\dfrac{1}{1+u}- \dfrac{1}{\mathrm{e}^{-\pi \omega}+ u}\right){\rm d} u \\
&=\dfrac{2\mathrm{e}^{\frac{-\pi \omega}{2}}}{\mathrm{e}^{-\pi \omega}-1} \log \left(\dfrac{1+u}{\mathrm{e}^{-\pi \omega}+ u}\right) \Big|_0^\infty 
= -\dfrac{2\mathrm{e}^{\frac{-\pi \omega}{2}}}{\mathrm{e}^{-\pi \omega}-1} \pi \omega \\
&=\pi \omega \text{csch}\left(\frac{\pi \omega}{2}\right).
\end{align*}
Therefore, $$\mathcal{F}\text{sech}^2 (\omega)=I(\omega)=\pi \omega
\text{csch}\left(\frac{\pi \omega}{2}\right)$$ and we conclude that the following
holds:
\begin{theorem}The univariate generalised Fourier transform of the
  hyperbolic tangent is
$$\mathcal{F}\bigl(\tanh x\bigr)(\omega)=\dfrac{1}{\mathrm{i}\omega}\mathcal{F}   \text{sech}^2 (\omega)=
\dfrac{1}{\mathrm{i}\omega}\pi \omega \text{csch}\,\left(\dfrac{\pi\omega}{2}\right) 
=-\mathrm{i}\pi\,\text{csch}\,\left(\dfrac{\pi\omega}{2}\right).$$
\end{theorem}

\subsection*{A.2. Generalised Fourier transforms of some classical radial basis functions}

Let ${{}}{x},{y}\in\mathbb{R}^n$ and $s=\| {y} \|$.

\begin{theorem*}{\cite[Th. 7.31]{jones1982}}
If $\beta\neq 2k$ and $\beta\neq -n-2k$ ($k=0,1,2,\ldots$) then the
$n$-dimensional generalised Fourier transform of
$r^\beta$ with $r=\|x\|$
is $$\frac{\Gamma\left(\frac{1}{2}\beta+\frac{1}{2}n\right)}{\Gamma\left(-\frac{1}{2}\beta\right)}
2^{\beta+n}\pi^{n/2}s^{-\beta-n},$$
where $s=\|y\|$ and $y\in\RR^n\setminus\{0\}$.
\end{theorem*}

\begin{theorem*}{\cite[Th. 7.32]{jones1982}}
For $k=0,1,\ldots$, the $n$-dimensional generalised Fourier transform of
$r^{2k}$ with $r=\|x\|$ is
$$ (2\pi)^n (-1)^k (\partial_1^2+\cdots+\partial_n^2)^k \delta
({{}}{y}),$$
where $y\in\RR^n\setminus\{0\}$.
\end{theorem*}

\begin{theorem*}{\cite[Th. 7.33]{jones1982}}
For $k=0,1,\ldots$, the $n$-dimensional generalised Fourier transform of
$r^{-n-2k}$ with $r=\|x\|$
is $$ \frac{(-1)^k\pi^{n/2}s^{2k}}{\Gamma\left(\frac{1}{2}n+k\right)
  \,k! \, 2^{2k-1}}\left\{ \frac{1}{2}\Psi \left(\frac{1}{2}n+k-1
\right)+\frac{1}{2}\Psi(k)+\log 2-\log s \right\},$$
where $s=\|y\|$ and $y\in\RR^n\setminus\{0\}$.
\end{theorem*}

\begin{theorem*}{\cite[Th. 7.34]{jones1982}}
  If $\beta\neq 2k$ and $\beta\neq -n-2k$ ($k=0,1,2,\ldots$) then
  the $n$-dimensional generalised Fourier transform of
  $r^\beta \log r $ 

  with $r=\|x\|$ is
  $${}\frac{\Gamma\left(\frac{1}{2}\beta+\frac{1}{2}n\right)}{\Gamma\left(-\frac{1}{2}\beta\right)}
2^{\beta+n}\pi^{n/2}s^{-\beta-n}\left\{ \frac{1}{2}\Psi
\left(\frac{1}{2}\beta+\frac{1}{2}n-1 \right)+\frac{1}{2}\Psi
\left(-\frac{1}{2}\beta-1 \right)-\log \frac{1}{2} s \right\},$$
where $s=\|y\|$ and $y\in\RR^n\setminus\{0\}$. $\Psi$ is the
Digamma-function, \cite{grad}.
\end{theorem*}

\begin{theorem}{\cite[Th. 7.34]{jones1982}}
For $k=0,1,\ldots$, the $n$-dimensional generalised Fourier transform of

$r^{2k} \log r$ 
is
$$\Gamma\left( \frac{1}{2}n+k\right) k! 2^{n+2k-1} \pi^{n/2} (-1)^{k+1} s^{-n-2k}
+
$$
$$(-1)^k(2\pi)^n\left\{\frac{1}{2}\Psi \left(\frac{1}{2}n+k-1 \right)+\frac{1}{2}\Psi (k)+\log 2\right\} (\partial_1^2+\cdots+\partial_n^2)^k \delta ({{}}{y}),
$$
where $s=\|y\|$ and $y\in\RR^n\setminus\{0\}$.

\end{theorem}

\begin{theorem}{\cite[Th. 7.35]{jones1982}}
For $k=0,1,\ldots$, the $n$-dimensional generalised Fourier transform of
$r^{-n-2k} \log r $ 
is
$$\frac{(-1)^k\pi^{n/2}s^{2k}}{\Gamma(\frac{1}{2}n+k) \,k! \,
  2^{2k}}\left[\left\{ \frac{1}{2}\Psi \left(\frac{1}{2}n+k-1
  \right)+\frac{1}{2}\Psi(k)-\log \frac{1}{2} s \right\}^2\right.$$
  $$
-\left.\frac{1}{4}\left\{\Psi ' \left(\frac{1}{2}n+k-1\right)+\Psi '
(k)\right\}+\frac{\pi^2}{12}\right],
$$
where $s=\|y\|$ and $y\in\RR^n\setminus\{0\}$.
\end{theorem}

\bibliographystyle{plain} 
\bibliography{arx}
\end{document}